\documentclass[12pt]{article}
\usepackage{color}
\usepackage{hyperref}

\usepackage{amsmath,amssymb,amsbsy}
\newcommand{\beq}{\begin{eqnarray}}
\newcommand{\eeq}{\end{eqnarray}}
\newcommand{\bei}{\begin{itemize}}
\newcommand{\eei}{\end{itemize}}
\newcommand{\bee}{\begin{enumerate}}
\newcommand{\eee}{\end{enumerate}}
\newcommand{\beqe}{\begin{eqnarray*}}
\newcommand{\eeqe}{\end{eqnarray*}}

\newtheorem{prop}{Proposition}
\newtheorem{lemma}{Lemma}

\newtheorem{theo}{Theorem}

\newcommand{\pa}[1]{\left({#1}\right)}
\newcommand{\cro}[1]{\left[{#1}\right]}

\newcommand{\ac}[1]{\left\{{#1}\right\}}

\newcommand{\R}{\mathbb{R}}

\newcommand{\N}{\mathbb{N}_*}

\newcommand{\ri}{\rightarrow}

\begin{document}

\renewcommand{\theequation}{\thesection.\arabic{equation}}
\title{Signal detection for inverse problems in a multidimensional framework\thanks{This work was realized during the invitation of Yuri 
Ingster at the INSA of Toulouse in June and July 
2012. He disappeared a few weeks later. We are very honored to have had the opportunity to work with him.}}
\author{ \framebox{Y. Ingster}, B. Laurent, C. Marteau }

\maketitle
\begin{abstract}
This paper is devoted to multi-dimensional inverse problems. In this setting, we address a goodness-of-fit testing problem. 
We investigate the separation rates associated to different kinds of smoothness assumptions
 and different degrees of ill-posedness. 
\end{abstract}

\section{Introduction}
This paper is concerned with an inverse problem model. More formally, given $\mathcal{H},\mathcal{K}$ two Hilbert spaces,
 our aim is to provide some inference on 
a function $f\in \mathcal{H}$ of interest from an observation $Y$ satisfying
\begin{equation}
Y= Kf + \epsilon \xi, 
\label{eq:model}
\end{equation}
where $K: \mathcal{H} \rightarrow  \mathcal{K}$ denotes a compact operator, $\epsilon$ a positive noise level and $\xi$ a Gaussian white noise. The model (\ref{eq:model}) means in fact that we can observe
\begin{equation}
\langle Y,g \rangle = \langle Kf, g \rangle + \epsilon \langle \xi, g \rangle, \ \forall g\in \mathcal{K},
\label{eq:model_angle}
\end{equation}
where for all $g,g_1,g_2 \in \mathcal{K}$, $\langle \xi, g \rangle \sim \mathcal{N}(0,\| g \|^2)$ and $\mathbb{E} \cro{\langle \xi, g_1 \rangle\langle \xi, g_2 \rangle } = \langle g_1,g_2 \rangle$. Such a model has been widely studied in the literature. In particular, the estimation issue has retained a large amount of attention. We mention for instance \cite{EHN_1996} for an extended review of existing methods and problematics in a numerical setting (i.e. when $\xi$ is deterministic and $\|\xi \| \leq 1$). In a statistical setting, we mention \cite{MP_2001}, \cite{Bissantz_2007} or \cite{Cavalier_book} for a recent survey on this topic.\\

In this paper, our aim is not to provide an estimator of $f$ and to study the related performances,
 but rather to  assess 
the goodness-of-fit of $f$ with respect to a benchmark signal $f_0$. More formally, our aim is to test
\begin{equation}
H_0: f=f_0, \ \mathrm{against} \ H_1: f\in \mathcal{F},
\label{eq:pb_test}
\end{equation}
where $\mathcal{F} \subset \mathcal{H}$ is such that $f_0 \not \in \mathcal{F}$. The description of the subset $\mathcal{F}$ is
 of first interest since it characterizes the detectable functions. This set is not allowed to be too rich since in this case it will
 be impossible to separate both hypotheses with prescribed errors. In the same time, it should not contain signals that are too close of 
$f_0$. Hence, most of the alternatives
in the literature are of the form
$$ \mathcal{F} := {\Theta}[r] = \left\lbrace f\in {\Theta}, \ \mathrm{s.t.} \ \| f- f_0\| \geq r \right\rbrace,$$
where ${\Theta}$ denotes a functional space  while $r$ describes the amount of signal available in the observations. Typically, 
  ${\Theta}$ will be a set of functions satisfying some smoothness constraint. 
Given ${\Theta}$ and a prescribed level $\alpha$ for the first kind error, the main challenge in this setting is to describe possible
 values of $r$ for which the second kind error can be (asymptotically) controlled. 
Without loss of generality, we assume in the following that $f_0=0$.  \\

In the direct setting (i.e. when $K$ denotes the identity operator), we mention the seminal series of paper \cite{IngsterI}, \cite{IngsterII}, 
\cite{IngsterIII} where various situations are considered. We also refer to \cite{Baraud} for a non-asymptotic study and \cite{BHL_2003} where adaptation with respect to  the smoothness of the alternative is discussed. In the inverse problem framework, few papers were interested by the testing issue. We mention \cite{Butucea} in a slightly different framework (error in variable model), or more recently \cite{LLM_2011}, \cite{LLM_2012}, \cite{ISS_2012} and \cite{MM_2013}  in our context.  

All these studies were concerned with uni-dimensional problem, where typically $\mathcal{H} = L^2(\mathbb{R})$. Up to our 
knowledge, the only study in the multi-dimensional case for inverse problem is provided in \cite{ISS_Radon}. Nevertheless, the 
framework is restricted to the Radon transform operator in dimension 2. In this paper, our aim is to investigate
 the behavior of separation rates in a d-dimensional setting, under general assumptions on the operator.\\

Our paper is organized as follows. In Section \ref{s:model}, we provide a precise description of the considered inverse problem model. 
We establish our main theorem, which will allow to determine the rates in the different considered cases. Section \ref{s:rates} is devoted 
to a presentation of the  minimax separation rates corresponding to various situations (mildly and severely ill-posed problems), 
regular and super-smooth alternatives. All the proofs are gathered in 
Section \ref{s:proof1}.

\section{A general methodology for testing in an inverse problem framework}
\label{s:model}

In this section, we will describe more precisely our model and the different related assumptions. Then we establish our main result which is at the core of our analysis.

\subsection{Multi-dimensional inverse problems}
The inverse model (\ref{eq:model}) has been widely investigated in the literature, in particular for an estimation purpose. Such a model 
arises in various practical situations and precise analysis are of first interest. Below, we describe two particular examples of operator and related applications.\\

\noindent
\textbf{The Radon transform.} The Radon transform is one of most popular operator. It is often involved in medical imaging. Given $H$ the unit disk on $\mathbb{R}^2$, the aim is to provide inference on the spatially varying density $f$ of a cross section $\mathcal{D} \subset H$ of an human body. This inference is provided via observation obtained by non-destructive imagery, namely X-ray tomography. In such a case, given a function $f \in L^2(H)$, one measure $Rf(u,\varphi)$ which corresponds to the decay of the intensity of the X-ray in the direction $\varphi$ when the receiver is at a distance $u$. Typically, we have
$$ Rf(u,\varphi) =\frac{\pi}{2\sqrt{1-u^2}} \int_{-\sqrt{1-u^2}}^{\sqrt{1-u^2}} f(u\cos(\varphi)-t\sin(\varphi),u\sin(\varphi)+t\cos(\varphi))dt,$$
for all $(u,\varphi) \in [0,1]\times[0,2\pi)$. For more details on such an operator, we mention for instance \cite{EHN_1996}. Testing issues are discussed in \cite{ISS_2012}. \\

\noindent
\textbf{Convolution operators.} The convolution operator $K$ is defined as follows
\begin{eqnarray*}
K :& & L^2(\mathbb{R}^d) \rightarrow L^2(\mathbb{R}^d)\\
& & f \mapsto Kf = \int_{\mathbb{R}^d} g(y) f(.-y)dy,
\end{eqnarray*}
where $g \in L^1(\mathbb{R}^d)$ denotes a convolution kernel. The behavior of this kernel (in particular the spectral behavior) characterizes the difficulty of the related deconvolution problem. Up to some conditions on this kernel $g$, the related convolution operator appears to a be a compact operator. For more details on the associated problem, we mention \cite{EHN_1996}, \cite{Cavalier_book} or \cite{Butucea} for the testing issue in a slightly different 
setting. \\

%\noindent
%\textbf{Differentiation as an inverse problem.}  .\\

The aim of this paper is not to concentrate on particular operators, but rather to describe a general scheme leading to satisfying testing procedures and related separation rates. Hence, in the following, we will only assume that we deal with the model (\ref{eq:model}) where $K$ denotes a compact operator. 

In order to determine whether $f=0$ or not (see \ref{eq:pb_test}), the main underlying idea is to construct an estimator of $\|f\|^2$ 
(or equivalently $\|Kf\|^2$, see discussion below) 
from $Y$ and 
then to take a decision based on this estimator. Indeed, if the corresponding norm is large enough  with respect to  some prescribed 
threshold, it seems reasonable to reject $H_0$. \\

In order to provide an  estimator of this norm, we will express our problem in terms of the singular value decomposition
 (S.V.D.) of the operator $K$. This decomposition is expressed through a sequence $(b_l^2,\phi_l,\varphi_l)_{l\in \N^d}$ where the
$b_l^2$ correspond to the eigenvalues of $K^*K$ while the $\phi_l$ correspond to the eigenfunctions (which are assumed to denote a basis of 
$\mathcal{H}$). In particular, for all $l\in \N^d$, we get
$$ \left\lbrace
\begin{array}{l}
K \phi_l = b_l \varphi_l, \\
K^* \varphi_l = b_l \phi_l,
\end{array} \right.$$
where $K^*$ denotes the adjoint operator of $K$. Using the previous relationships and (\ref{eq:model_angle}), we obtain observations as a sequence $(y_l)_{l\in \N^d}$ with
\begin{equation}
y_l = b_l \theta_l + \epsilon \xi_l, \ l \in \N^d,
\label{eq:sequence_space_model}
\end{equation}
where for all $l\in \N^d$, $\theta_l = \langle f,\phi_l \rangle$ and the $\xi_l$ are i.i.d. standard Gaussian random variables. 
The model (\ref{eq:sequence_space_model}) is usually called the sequence space model. The main advantage of such a representation is that 
it provides a simple way to investigate the different interactions between properties of the function of interest and the related behavior of the operator. Following \cite{ISS_2012} or \cite{LLM_2011}, we will only deal in the following with the model (\ref{eq:sequence_space_model}). 

Please note that it is also possible to adopt another point of view where one may use general regularization approaches instead of projecting 
observations onto particular bases (here the S.V.D. basis). In  the testing issue, this methodology has been investigated in a
uni-dimensional inverse problem framework in \cite{MM_2013}. Nevertheless, such an approach appears to be more complicated in the
 multidimensional case and hence will not be considered in this paper.

\subsection{Minimax separation rates}
Let $$ {\Theta}(r_\epsilon) =  \left\lbrace f\in {\Theta}, \ \mathrm{s.t.} \ \| f \| \geq r_\epsilon \right\rbrace,$$
where $r_\epsilon$ denotes some radius, which is allowed to depend on $\epsilon$. The set $\Theta$ is a given functional 
space (which will be made precise later on). Our aim in this paper is to test 
\begin{equation}
H_0: f=0, \ \mathrm{against} \ H_{1,r_\epsilon}: f\in {\Theta}(r_\epsilon).
\label{eq:intro_test}
\end{equation}
For this purpose, we have to propose a test function $\psi$ which is a  measurable
 function of the data with values in $\lbrace 0,1 \rbrace$. By convention, we reject $H_0$ if $\psi=1$, and do not reject
 otherwise. Given a prescribed level $\alpha \in ]0,1[$ for the first kind 
error, we will deal all along the paper with level-$\alpha$ tests, namely test functions $ \psi_\alpha$ satisfying
$$ \mathbb{P}_{H_0} (\psi_\alpha =1) \leq \alpha + o(1), \ \mathrm{as} \ \epsilon \rightarrow 0.$$
Then, we associate to each  level-$\alpha$ test function $\psi_\alpha$ its maximal type II  error probability over the set 
${\Theta}(r_\epsilon)$  defined as
$$ \beta_{\epsilon}({\Theta}(r_\epsilon),\psi_\alpha)= \sup_{f\in {\Theta}(r_\epsilon)} \mathbb{P}_f( \psi_\alpha =0). $$ 
The minimax type II error probability over the set ${\Theta}(r_\epsilon)$ is defined as
$$ \beta_{\epsilon}({\Theta}(r_\epsilon),\alpha) =\inf_{\psi_\alpha}  \beta_{\epsilon}({\Theta}(r_\epsilon),\psi_\alpha),$$ 
where the infimum is taken over all possible  level-$\alpha$ test functions. \\

Given a test function $\psi_\alpha$ and a radius $r_{\epsilon}$,  $ \beta_{\epsilon}({\Theta}(r_\epsilon),\psi_\alpha)$ hence provides 
a measure of the performances of the test. 
Typically, if this term tends to $1-\alpha $ as $\epsilon \rightarrow 0$, the test $\psi_\alpha$ does not separate both hypotheses. 
On the other hand, given any $\beta \in ]0,1[$, if one can find a test function $\psi_{\alpha} $ such that
 $ \beta_{\epsilon}({\Theta}(r_\epsilon),\psi_{\alpha})\leq \beta + o(1)$ as $\epsilon \rightarrow 0$, then 
the  hypotheses $H_0$ and $  H_{1,r_\epsilon} $ can be separated with prescribed levels. 

For a given  radius   $r_{\epsilon}$, the term $ \beta_{\epsilon}({\Theta}(r_\epsilon),\alpha)$ represents the lowest achievable 
level for the type II
error probability. Clearly, this quantity is non-increasing with respect to $r_{\epsilon}$. We are therefore interested in the minimal 
$r_{\epsilon}$ for which $  \beta_{\epsilon}({\Theta}(r_\epsilon),\alpha) \rightarrow 0 $ as $\epsilon \rightarrow 0$. 
In particular, a sequence $r_{\epsilon}^\star$ is called a minimax separation rate over the set $\Theta$ if 
\begin{eqnarray}
&&  \beta_{\epsilon}({\Theta}(r_\epsilon),\alpha) \rightarrow 1-\alpha \hspace{0.5cm} \mbox{ if } r_{\epsilon}/r_{\epsilon}^\star  \rightarrow 0 \label{minimax1} \\
&&  \beta_{\epsilon}({\Theta}(r_\epsilon),\alpha) \rightarrow 0 \qquad \hspace{0.5cm} \mbox{ if } r_{\epsilon}/r_{\epsilon}^\star  \rightarrow \infty \label{minimax2} 
\end{eqnarray}

Hence, from an asymptotic point of view, the term $r_\epsilon^\star$ denotes the order of the delimiting separation radius under which 
testing is impossible (i.e. for which the type II error probability will be close to $1-\alpha$). 

One of the main issue in the testing theory is then to describe as precisely as possible the (asymptotic) value for the separation rate 
$r_\epsilon^\star$ following the smoothness constraints expressed on the signal. In an inverse problem framework, one want also to take 
into account the behavior of the operator when describing $r_\epsilon^\star$. In the following, we introduce our different constraints on 
both the signal and the operator. Then, we will establish a general result that will allow us to investigate the asymptotic of 
$r_\epsilon^\star$ in different settings. \\

Following classical results in an inverse problem framework (see for instance \cite{EHN_1996} or \cite{Cavalier_book}), the difficulty of the problem will be characterized by the behavior of the sequence $(b_l)_{l\in \N^d}$. Indeed, starting from the model (\ref{eq:sequence_space_model}), it appears that it will be difficult
 to retrieve informations on  $\theta_l$ if the corresponding coefficient $b_l$ is close to $0$. Basically, two different regimes are 
considered in the literature. \\
\\
\textbf{Mildly ill-posed problems} \textit{There exists $0<c_0\leq c_1$ such that the sequence $(b_l)_{l\in \N^d}$ satisfies
$$ c_0  \prod_{j=1}^d |l_j|^{-t_j} \leq |b_l| \leq c_1 \prod_{j=1}^d |l_j|^{-t_j}, \ \forall l=(l_1,\ldots, l_d) \in \N^d,$$ 
for some sequence $t=(t_1,\dots,t_d)\in \R_+^d$.}\\

Mildly ill-posed inverse problems correspond to the most favorable cases in the sense that the sequence $(b_l)_{l\in \N^d}$ does not decrease too fast. On the other hand, severely ill-posed problems are more difficult to handle.\\
\\
\textbf{Severely ill-posed problems} \textit{ There exists $0<c_0\leq c_1$ such that the sequence $(b_l)_{l\in \N^d}$ satisfies  
$$  c_0 \prod_{j=1}^d e^{-t_j l_j} \leq  |b_l| \leq c_1 \prod_{j=1}^d e^{-t_j l_j}, \ \forall l= (l_1,\ldots, l_d) \in \N^d,$$ 
for some sequence $(t_1,\dots,t_d)\in \R_+^d $.}\\

The previous conditions characterize most of the inverse problems encountered in the literature. Now, we have to introduce smoothness assumptions on our target $f$. In the following, given a  sequence $a=(a_l)_{l\in \N^d}$ and a positive radius $R$, we will use the corresponding ellipsoids 
${\Theta}$ defined as
$$ \Theta:={\Theta}_{a,R} = \left\lbrace \nu \in l^2(\N^d), \ \sum_{l\in \N^d} a_l^2 \nu_l^2 \leq R^2 \right\rbrace.$$
In particular, two different kinds of behavior for the sequence $a=(a_l)_l$ will be considered. Given positive numbers $(s_1,\ldots , s_d)$, 
we will alternatively deal with sequences $a$ satisfying
\begin{equation}\label{aSobolev}
 a_l^2 = \sum_{j=1}^d l_j^{2s_j}, \quad \ \mathrm{or} \quad \ a_l^2= \sum_{j=1}^d e^{2s_j l_j} \ \quad \forall l\in \N^d,
 \end{equation} 
and sequences $a$ satisfying
\begin{equation}\label{aTensor}
 a_l^2 = \prod_{j=1}^d l_j^{2s_j}, \  \quad \ \mathrm{or} \quad \ a_l^2= \prod_{j=1}^d e^{2s_j l_j} \ \quad \forall l\in \N^d.
 \end{equation}
 In the first case (\ref{aSobolev}), the ellipsoid $ {\Theta}_{a,R}$  corresponds to a \textit{Sobolev} space, while in the second 
case (\ref{aTensor}), we deal with a so-called \textit{tensor product} space.

\subsection{Main result}

In this section, we provide a general result that will allow to determine minimax separation rates in various situations. 
As formalized in (\ref{eq:intro_test}), our aim is to decide whether there is signal in our observations or not. A natural way to investigate
 this problem is to construct a criterion that will measure the amount of signal available in the observations. In this paper, our test will 
be based  on a estimation of the norm of $Kf$. Indeed, assertions $"\|f\| =0"$, $"f=0"$ and $"Kf=0"$ are equivalent since our operator is 
 injective (we mention \cite{LLM_2011} for an extended discussion). \\

First, we can estimate $\|Kf\|^2$ using the statistics 
\begin{equation}
T = \sum_{k\in \N^d} \omega_k (y_k^2 - \epsilon^2),
\label{eq:est_filtre}
\end{equation}
where $\omega=(\omega_k)_{k\in \N^d}$ denotes a filter, i.e. a real sequence with values in $[0,1]$. Several kinds of filters are 
available in the literature for a testing purpose. We can mention for instance \cite{Baraud}, \cite{LLM_2011} or \cite{LLM_2012} where 
projection filters designed as $\omega_k = \mathbf{1}_{\lbrace k_1\leq N_1,..,k_d\leq N_d\rbrace}$   are proposed or \cite{MM_2013} where 
properties of  Tikhonov filters (and more general regularization schemes) are investigated.  From now on, following seminal investigations 
proposed in \cite{IngsterI}-\cite{IngsterIII}, and more recent results in an inverse problem setting, we will consider  filters  designed as
\begin{equation}
\omega_k = \frac{b_k^2\tilde \theta_k^2}{\sqrt{2\sum_{k\in \N^d} b_k^4\tilde\theta_k^4}} \ \forall k\in \N^d,
\label{eq:filtre_ingster}
\end{equation}
where for a given $r_\epsilon$ the sequence $(\tilde \theta_k)_k$ is the solution of the extremal problem 
\begin{equation}
 u_\epsilon^2 (r_\epsilon) = \frac{1}{2\epsilon^4} \sum_{k\in \N^d} b_k^4 \tilde \theta_k^4 = \frac{1}{2\epsilon^4}
 \inf_{\theta \in \Theta_{a,R}(r_\epsilon)} \sum_{k\in \N^d}   b_k^4 \theta_k^4.
\label{eq:extremal_problem}
\end{equation}
In general, a decision rule is based on the following principle. If the norm of the estimator of $Kf$ is larger than a given threshold, this means that there is probably signal in
 the observations and we will reject $H_0$. On the other hand, if the norm of this estimator is not large enough, we are observing only noise with a
 high probability and we do not reject $H_0$. Using this principle and (\ref{eq:est_filtre})-(\ref{eq:filtre_ingster}), we obtain the testing procedure $\Psi_{\epsilon,t}$ defined as 
\begin{equation}
\Psi_{\epsilon,t} = \mathbf{1}_{\lbrace T> t \rbrace}, \ \mathrm{where} \ T= \sum_{k\in \N^d} \omega_k (y_k^2 - \epsilon^2),
\label{eq:testing_procedure}
\end{equation}
for some threshold $t$. In the following, for a given $\alpha \in ]0,1[ $, we set $t=H^{(\alpha)}$, where $H^{(\alpha)}$
denotes the $1-\alpha$ quantile of the standard Gaussian distribution. We will prove that the corresponding test is asymptotically of 
level $\alpha$, which means that
$$ \mathbb{P}_{H_0} (\Psi_{\epsilon,H^{(\alpha)}} =1)\leq  \alpha + o(1) \mbox{ as } \epsilon \rightarrow 0 .$$  

The following theorem emphasizes the performances of the test $\Psi_{\epsilon,t}$. We also provide lower bounds that asses the optimality of this testing procedure. 

\begin{theo}
\label{thm:main}
Consider the testing problem introduced in (\ref{eq:intro_test}) and the testing procedure $\Psi_{\epsilon,t}$ defined 
in (\ref{eq:testing_procedure}).  Then, given $\alpha \in ]0,1[ $,
\begin{enumerate}
\item 
\begin{enumerate}
\item If $u_\epsilon (r_\epsilon) \rightarrow 0$, then $  \beta_{\epsilon}({\Theta_{a,R}}(r_\epsilon),\alpha) \rightarrow 1-\alpha$ as
 $\epsilon \rightarrow 0$. In this case, minimax testing is
 impossible. 
%If $u_\epsilon = \mathcal{O}(1)$, then $\lim \inf \beta_\epsilon(r_\epsilon,\alpha)>0$ and minimax consistent testing is impossible. 
\item If $u_\epsilon(r_\epsilon) = \mathcal{O}(1)$ and $\omega_0= o(1)$ as $\epsilon \rightarrow 0$, then 
 the test $\Psi_{\epsilon,H^{(\alpha)}}$ is a level-$\alpha$ test and is asymptotically 
minimax, i.e.
%$$ \alpha_\epsilon (\Phi_{\epsilon,H^{(\alpha)}}) \leq \alpha + o(1),$$
%$$ \beta_\epsilon(\Theta(r_\epsilon), \Phi_{\epsilon,u_\epsilon/2}) = \beta(r_\epsilon,\alpha)+o(1),$$
$$ \beta_\epsilon(\Theta_{a,R}(r_\epsilon),\Psi_{\epsilon,H^{(\alpha)}}) = \beta_{\epsilon}({\Theta_{a,R}}(r_\epsilon),\alpha) + o(1), \quad \mathrm{as} \ \epsilon \rightarrow 0.$$
Moreover,  we obtain the sharp asymptotics
$$ \beta({\Theta_{a,R}}(r_\epsilon),\alpha) = \Phi( H^{(\alpha)} - u_\epsilon (r_\epsilon))+ o(1), \quad \mathrm{as} \ \epsilon \rightarrow 0.$$
%$$ \gamma_{\epsilon}({\Theta}(r_\epsilon)) = 2 \Phi(-u_\epsilon/2)+ o(1).$$
\end{enumerate}
\item If $u_\epsilon(r_\epsilon)\rightarrow + \infty$, then the family of tests (\ref{eq:testing_procedure}) with $t=cu_\epsilon(r_\epsilon) $ for some
 $c\in ]0,1[$ are asymptotically consistent, i.e. $\beta(\Theta_{a,R}(r_\epsilon),\Psi_{\epsilon,t}) \rightarrow 0$ as $\epsilon \rightarrow 0$.
\end{enumerate}
\end{theo}
For the sake of brevity, we will not provide a complete proof since it follows the same lines than previous one established in the direct case (see \cite{IngsterI}-\cite{IngsterIII}) or in a uni-dimensional inverse setting by \cite{ISS_2012}.
Nevertheless, we will 
provide the main underlying ideas in Section \ref{s:preuve_main} below. \\

The main consequence of Theorem \ref{thm:main} is that the investigation of the minimax separation rate associated to the testing problem (\ref{eq:intro_test}) reduces to the study of the extremal problem (\ref{eq:extremal_problem}). In particular, the behavior of the term $u_\epsilon$ will provide meaningful informations on the difficulty of the testing problem. 

In Section \ref{s:rates} below, we investigate the behavior of this extremal problem and related separation rates for different kinds of 
smoothness assumptions and different degrees of ill-posedness.

\section{Separations rates}
\label{s:rates}

Following Theorem \ref{thm:main}, it appears that the extremal problem (\ref{eq:extremal_problem}) is of first importance for 
a precise understanding of minimax separation rates in this setting. In particular, we have to make explicit the terms $r_\epsilon$ for which
 $u_\epsilon = \mathcal{O}(1)$ as $\epsilon \rightarrow 0$. \\

To this end, we can remark in a first time that the solution $\theta^\star$ of (\ref{eq:extremal_problem}) is of the form
$$ (\theta_k^\star)^2 = z_0^2 b_k^{-4} (1-Aa_k^2)_+, \ \forall k\in \mathbb{N}^d,$$
where the terms $z_0=z_{0,\epsilon}$ and $A=A_\epsilon$ are determined by the equations
$$ \left\lbrace
\begin{array}{l}
r_\epsilon^2 = z_0^2 J_1,\\
1 = z_0^2 A^{-1} J_2,
\end{array}
\right.$$
with
\begin{eqnarray*}
J_1 & = & \sum_{k\in \N^d} b_k^{-4} (1-Aa_k^2)_+,\\
\mathrm{and} \quad J_2 & = & A\sum_{k\in \N^d} a_k^2 b_k^{-4} (1-Aa_k^2)_+.
\end{eqnarray*}
In particular, 
$$ u_\epsilon^2(r_\epsilon) = \epsilon^{-4} z_0^4 J_0/2, \ \mathrm{where} \ J_0 = J_1-J_2 = \sum_{k\in \N^d} b_k^{-4} (1-Aa_k^2)_+^2.$$
Using this methodology, we get separation rates over both \textit{Tensor product} and \textit{Sobolev} spaces, when considering alternatively \textit{mildly} and \textit{severely} ill-posed problems. These rates are summarized in Tables \ref{tab:tensor} and \ref{tab:sobolev}. The formal results, and related proofs are made explicit in the sections below.

In the following, we use the notation $ v_{\epsilon} \sim w_{\epsilon}$ as $\epsilon \rightarrow 0$ if there exist two constants $c_0$ and $c_1$ such that
$\forall \epsilon >0$, $c_0 w_{\epsilon} \leq v_{\epsilon} \leq c_1 w_{\epsilon}$.

\begin{table}[t]
\begin{center}
\begin{tabular}{lcc}
\hline
                    & \textbf{Mildly ill-posed}               & \textbf{Severely ill-posed} \\
                    &    $|b_l| = \prod_{j=1}^d |l_j|^{-t_j}$                    &    $|b_l| = \prod_{j=1}^d e^{-t_j l_j}$ \vspace{0.1cm} \\ 
\hline
$a_l^2 = \prod_{j=1}^d |l_j|^{2s_j}$   &  $\epsilon^{\frac{4}{4+c_{1}}}$ \hspace{3.3cm} \textsf{(i)}   & $\left( \frac{1}{4t_1} \ln \left( \frac{1}{\epsilon^4} \right)\right)^{-2s}$ \quad \textsf{(iii)} \vspace{0.1cm} \\
\hline
$a_l^2 = \prod_{j=1}^d e^{2s_j l_j}$ &  $\epsilon^2 \left( \ln(\epsilon^{-1}) \right)^{\sum_{j=1}^d(2t_j +1/2)}$ \quad \textsf{(ii)}  & $\epsilon^{\frac{2}{1+t_1/s_1}}$ \hspace{1.8cm} \textsf{(iv)} \vspace{0.1cm} \\
\hline
\end{tabular}
\end{center}
\caption{\textit{Minimax separation rates for Tensor product spaces. In the case \textsf{(i)}, $c_1=\max_{j=1\dots d}(1+4t_j)/s_j$; in the case \textsf{(ii)}, $s_j=s$ for all $j\in\lbrace 1,\dots, d\rbrace$; in the case \textsf{(iii)}, $s_j=s$ and $t_1> t_j$
and $(t_1,s_1)$ for all $j\in\lbrace 1,\dots, d\rbrace$; in the case \textsf{(iv)},  $t_1/s_1 > t_j/s_j$ for all $j > 1 $. }}
\label{tab:tensor}
\end{table}

\begin{table}[t]
\begin{center}
\vspace{1.5cm}
\begin{tabular}{lcc}
\hline
                    & \textbf{Mildly ill-posed}               & \textbf{Severely ill-posed} \\
                    &    $|b_l| = \prod_{j=1}^d |l_j|^{-t_j}$                    &    $|b_l| = \prod_{j=1}^d e^{-t_j l_j}$ \vspace{0.1cm} \\ 
\hline
$a_l^2 = \sum_{j=1}^d |l_j|^{2s_j}$   &  $\epsilon^{{4}/\pa{4+\sum_{j=1}^d \frac{(1+4t_j)}{s_j}}}$   & $\left( \log \left(\frac{1}{\epsilon}\right) \right)^{-s}$ \vspace{0.1cm} \\
\hline
\end{tabular}
\end{center}
\caption{\textit{Minimax separation rates for Sobolev spaces. In the severely ill-posed case, we assume that $s_1=\dots=s_d=s$. }}
\label{tab:sobolev}
\end{table}

\subsection{Tensor product spaces}

\subsubsection{Mildly ill-posed problems with ordinary smooth functions}
In this section, we assume that 
\begin{equation}
b_l^{2} =\prod_{j=1}^d |l_j|^{-2t_j} \ \mathrm{and} \ a^2_{l}=\prod_{j=1}^d |l_1|^{2s_j}, \ \forall l\in \N^d.
\label{eq:cas1}
\end{equation}
In the following, we will deal with the sequence $(c_j)_{j=1\dots d}$ defined as
$$ c_j = \frac{1+t_j}{s_j}, \ \forall j\in \lbrace 1,\dots,d \rbrace.$$
The following proposition describes the minimax separation rate in this setting.

\begin{prop}
\label{prop:cas1}
Assume that both sequences $a$ and $b$ satisfy equation (\ref{eq:cas1}) and that $c_1,\ldots, c_d$ are strictly ordered :   
$c_1>\dots >c_d$. Then, we get that
\begin{itemize}
\item The minimax separation rate $r_\epsilon^\star$ satisfies
$$ r_\epsilon^\star \sim \epsilon^{\frac{4}{4+c_1}} \ \mathrm{as} \ \epsilon\rightarrow 0.$$
\item  If $ r_\epsilon = C  \epsilon^{\frac{4}{4+c_1}}$ for some positive constant $C$, then 
$u_{\epsilon}(r_\epsilon ) = \mathcal{O}(1)$ as $\epsilon \rightarrow 0$ and  the sharp asymptotics described in Theorem \ref{thm:main} hold true.
\end{itemize} 
\end{prop}
The proof of this result is postponed to Section \ref{s:preuve_cas1}.\\

Remark that in the particular case where $d=1$, we get the minimax rate
$$ r_\epsilon^\star = \epsilon^{\frac{4s}{4s+4t+1}},$$
which has been established for instance in \cite{ISS_2012} or \cite{LLM_2011}. Hence, Proposition \ref{prop:cas1} appears to be an extension of this classical $1$-dimensional case. Nevertheless, remark that the obtained rate appears to be quite unusual
 in such a setting: it is governed by the couple of parameters $(t_j,s_j)$ for which the associated term $c_j$ is maximal. 

In some sense, the rate corresponds to the worst $1$-dimensional rate in each direction. In particular, a large value for $c_j$ is more or less associated to high degree of ill-posedness with a small smoothness index.  

\subsubsection{Mildly ill-posed problems with supersmooth functions}
In this section, we assume that
\begin{equation}
b_l^2 = \prod_{j=1}^d l_j^{-2t_j} \ \mathrm{and} \ a_l^2 = \prod_{j=1}^d e^{2s l_j} \ \forall l\in \N^d.
\label{eq:cas2}
\end{equation}
In such a setting, we get the following result, whose proof is postponed to Section \ref{s:preuve_cas2}.

\begin{prop}
\label{prop:cas2}
Assume that both sequences $a$ and $b$ satisfy equation (\ref{eq:cas2}). Then, we get that
\begin{itemize}
\item The minimax separation rate $r_\epsilon^\star$ satisfies
$$ r_\epsilon^\star \sim \epsilon \left( \ln(\epsilon^{-1}) \right)^{\sum_{j=1}^d(t_j +1/4)} \ \mathrm{as} \ \epsilon\rightarrow 0.$$
\item  If $ r_\epsilon = C \epsilon \left( \ln(\epsilon^{-1}) \right)^{\sum_{j=1}^d(t_j +1/4)}  $ for some positive constant $C$, then 
$u_{\epsilon}(r_\epsilon ) = \mathcal{O}(1)$ as $\epsilon \rightarrow 0$ and  the sharp asymptotics described in Theorem \ref{thm:main} hold true.
\end{itemize} 
\end{prop}

In this setting, the minimax separation rate is close to the parametric separation rate, up to a logarithmic term. The power of this log term 
explicitly depends on the dimension, and of the degrees of ill-posedness in each direction. Once again, we recover the results obtained 
in \cite{ISS_2012} or \cite{LLM_2011} in the particular case where $d=1$.

\subsubsection{Severely ill-posed problems with super-smooth functions}
In this section, we deal with severely ill-posed problems and super-smooth functions. In particular 
\begin{equation} 
b_l = \prod_{j=1}^d e^{-t_j l_j} \ \mathrm{and} \ a_l = \prod_{j=1}^d e^{s_j l_j} \ \forall l\in \N^d.
\label{eq:cas3}
\end{equation}
The proof of the following proposition is provided in Section \ref{s:preuve_cas3}. 

\begin{prop}
\label{prop:cas3}
Assume that both sequences $a$ and $b$ satisfy equation (\ref{eq:cas3}) where
\begin{equation} 
\frac{t_1}{s_1} > \frac{t_j}{s_j} \ \forall j\in \lbrace 2,\dots,d \rbrace.
\label{eq:hyp_cas3}
\end{equation}
Then, the minimax separation rate $r_\epsilon^\star$ satisfies
$$ r_\epsilon^\star \sim \epsilon^{\frac{s_1}{s_1+t_1}} \ \mathrm{as} \ \epsilon\rightarrow 0.$$
\end{prop}
The assumption (\ref{eq:hyp_cas3}) appears to be necessary since it allows a sharp control of the terms $J_0,J_1$ and $J_2$ introduced above. It could be certainly removed, up to very technical algebra. In particular, we mention that the detection of super-smooth functions with severly ill-posed problems is already a very difficult problem in $1$-dimensional case. We also mention that in this particular setting, we do not obtain sharp asymptotics.

\subsubsection{Severely ill-posed problems with ordinary smooth functions}

In this section, we deal with severely ill-posed problems with ordinary isotropic smooth functions. In particular
\begin{equation}
b_l = \prod_{j=1}^d e^{-t_j l_j} \ \mathrm{and} \  a_l^2 = \prod_{j=1}^d l_j^{2s} \ \forall l\in \N^d.
\label{eq:cas4}
\end{equation}
The terms $(t_1,\dots,t_d)$ and $s$ denote positive known parameters, that characterize the problem. The proof of the following proposition is provided in Section \ref{s:preuve_cas4}.

\begin{prop}
\label{prop:cas4}
Assume that both sequences $a$ and $b$ satisfy equation (\ref{eq:cas4}) where
\begin{equation}
 t_1 > t_j \ \forall j\in \lbrace 2,\dots, d \rbrace.
\label{eq:cond_prop4}
\end{equation}
Then, the minimax separation rate $r_\epsilon^\star$ satisfies
$$r_\epsilon^\star \sim \left( \frac{1}{4t_1} \ln \left( \frac{1}{\epsilon^4} \right)\right)^{-2s} \ \mathrm{as} \ \epsilon\rightarrow 0.$$
\end{prop}

We work in an isotropic context, in the sense that the regularity is the same in all directions. Provided that property (\ref{eq:cond_prop4}) holds, we obtain the same rate that in an uni-dimensional framework. Then, the minimax separation rates is characterized by the direction for which the problem is the more difficult, i.e. the direction associated to the largest indice $t_j$. Once again, we do not get sharp rates in this setting.

\subsection{Sobolev Spaces}
In this part, we will consider Sobolev smoothness constraints.
 We will provide some illustrations of Theorem \ref{thm:main} in this context.
In both considered cases, we precise the behavior of the sequence $a$.

\subsubsection{Mildly ill-posed problems}
We start our study in this framework with mildly ill-posed inverse problems. In particular, we assume that
\begin{equation}
b^2_{l}= \prod_{j=1}^d l_1^{-2t_j} \ \mathrm{and} \ a^2_{l}= \sum_{j=1}^d l_j^{2s_j}, \ \forall l\in \N^d.
\label{eq:sobolev_cas1}
\end{equation}

\begin{prop}
\label{prop:sobolev_cas1}
Assume that both sequences $a$ and $b$ satisfy (\ref{eq:sobolev_cas1}). Then, we have that 
\begin{itemize}
\item The minimax separation rate satisfies
$$ r_{\epsilon}^* \sim \epsilon^{{4}/\pa{4+\sum_{j=1}^d \frac{(1+4t_j)}{s_j}}} \ \mathrm{as} \ \epsilon \rightarrow 0.$$
\item If $ r_{\epsilon} =C \epsilon^{{4}/\pa{4+\sum_{j=1}^d \frac{(1+4t_j)}{s_j}}} $ 
for some positive constant $C$, then  $u_{\epsilon}(r_\epsilon ) = \mathcal{O}(1)$ as $\epsilon \rightarrow 0$ and 
the sharp asymptotics described in Theorem \ref{thm:main} hold true.

\end{itemize}
\end{prop}
The proof of Proposition \ref{prop:sobolev_cas1} is postponed to the Section \ref{s:sobolev_cas1}. \\

In this context, some cases are or first interest. In a first time, we can remark that in dimension $1$, we recover the classical rate
$$ r_{\epsilon}^*=  \epsilon^{\frac{4s}{4s+4t+1}}.$$
We refer for more details to \cite{ISS_2012} or \cite{LLM_2012} (and to the previous section for a related discussion). The other interesting case corresponds to the homoscedastic framework in dimension $d$ where $s_j=s$ and $t_j=t$ for all $j\in \lbrace 1,\dots, d \rbrace$. In such a case, we get that
$$ r_{\epsilon}^*=  \epsilon^{\frac{4s}{4s+4dt+d}}.$$
In this situation, we are faced to the well-known curse of the dimension: minimax detection rates deteriorate as $d$ increases, even for direct problems (where $t=0$).

\subsubsection{Severly ill-posed problems}
We now turn to the investigation of severly ill-posed inverse problems. We will consider the case where
\begin{equation} \label{SIPP}
 b_l = \prod_{j=1}^d e^{-t_j l_j} , \ \mathrm{and} \ a_l = \left(\sum_{j=1}^d l_j \right)^s, \ \forall l\in \N^d.
 \end{equation}
In particular, we only deal with homoscedastic problems: the regularity is supposed to be the same in all considered directions. 
In this context, we get the following result, whose proof can be found in Section \ref{s:sobolev_cas2}.

\begin{prop}
\label{prop:sobolev_cas2}
Assume that both sequences $a$ and $b$ satisfy (\ref{SIPP}), where $t_1>t_j$ for all $j\in \lbrace 2,\dots, d \rbrace$.
 Then, we obtain the following results :
\begin{itemize}
 \item If $ r_\epsilon = \pa{C\log(1/\epsilon)}^{-s} $ with $C>1/t_1$, then $ u_\epsilon^2 (r_\epsilon)\rightarrow +\infty$ as $\epsilon \rightarrow
0$ and the detection is possible (see Theorem \ref{thm:main} for more details).
\item If $ r_\epsilon = \pa{C\log(1/\epsilon)}^{-s} $ with $C\leq 1/t_1$, then $ u_\epsilon^2 (r_\epsilon)\rightarrow 0 $ 
as $\epsilon \rightarrow 0$ and the detection is impossible. 
\end{itemize}
\end{prop}
Due to the difficulty of the problem, the minimax separation rate decreases very slowly and we do not obtain sharp asymptotics in this case.  
In order to get this result, we have assumed that the problem is in some sense dominated by one direction:
 $t_1>t_j$ for all $j\geq 2$. This assumption can be removed, up to a more technical algebra.

\section{Proofs}
\label{s:proof1}

\subsection{Proof of Theorem \ref{thm:main}}
\label{s:preuve_main}
The proof is decomposed in two different parts. In a first time, we establish a lower bound for the term $ \beta_{\epsilon}({\Theta}(r_\epsilon),\alpha)$ and we discuss the possible values of this quantity following the behavior of the extremal problem (\ref{eq:extremal_problem}). Then, we prove that the test $\Psi_{\epsilon,H^{(\alpha)}}$ achieves this lower bound when $u_\epsilon(r_\epsilon) = \mathcal{O}(1)$ as $\epsilon\rightarrow 0$.  \\

We first focus on the lower bound on $ \beta_{\epsilon}({\Theta}(r_\epsilon),\alpha)$. Let $\pi$ the prior on the set $\Theta(r_\epsilon)$ defined as 
\begin{equation}
\pi = \prod_{k\in \N^d} \pi_k, \ \mathrm{where} \ \pi_k=\frac{1}{2} (\delta_{-b_k\theta_k} + \delta_{b_k\theta_k}) \ \forall k\in \N^d,
\label{eq:measure}
\end{equation}
for some sequence $\theta \in \Theta(r_\epsilon)$ which will  be made precise later on. Denote by $\mathbb{P}_0$ (resp. $\mathbb{P}_\pi$) the measure associated to the observation vector $Y$ when the sequence $\theta$ is equal to $0$ (resp. follows the measure $\pi$). Then, following \cite{Baraud}, we get that
$$  \beta_{\epsilon}({\Theta}(r_\epsilon),\alpha) \geq 1- \alpha -\frac{1}{2} \left( \mathbb{E}_0[ L_{\pi}^2(Y)] -1 \right)^{1/2},$$
where $L_\pi(Y)$ denotes the likelihood ratio between the two measures $\mathbb{P}_0$ and $\mathbb{P}_\pi$. Thanks to (\ref{eq:measure}), we get that
$$\mathbb{E}_0[ L_{\pi}^2(Y)] = \prod_{k\in\N^d} \cosh (b_k^2\theta_k^2/\epsilon^2) \leq \exp \left( \frac{1}{2\epsilon^4} \sum_{k\in \N^d} b_k^{4}\theta_k^{4} \right) := 
\exp(u_\epsilon^2(r_\epsilon)),$$
provided $\theta$ is defined as the solution of the extremal problem (\ref{eq:extremal_problem}). In this context, we get clearly that
$$ u_\epsilon^2(r_\epsilon) \rightarrow 0 \Rightarrow \beta_{\epsilon}({\Theta}(r_\epsilon),\alpha) \rightarrow 1-\alpha \hspace{1cm} \mathrm{as} \ \ \epsilon \rightarrow 0,$$
which prove item 1.(a) of the Theorem. \\

Now, we turn to the case where $u_\epsilon(r_\epsilon) = \mathcal{O}(1)$ as $\epsilon \rightarrow 0$. In that case, our aim is to prove that 
$$ \beta_{\epsilon}({\Theta}(r_\epsilon),\alpha) =  \Phi(H^{(\alpha)}-u_\epsilon) + o(1) \ \mathrm{as} \ \epsilon \rightarrow 0.$$
To this end, assume that 
\begin{equation}
\ln (L_\pi(Y)) = -\frac{u_\epsilon^2(r_\epsilon)}{2} + u_\epsilon(r_\epsilon) Z_\epsilon + \delta_\epsilon,
\label{eq:rapport_vrais}
\end{equation}
where $Z_\epsilon \rightarrow Z\sim \mathcal{N}(0,1)$ and $\delta_\epsilon \rightarrow 0$ in $\mathbb{P}_0$-probability as $\epsilon \rightarrow 0$. It is well known that
$$  \beta_{\epsilon}({\Theta}(r_\epsilon),\alpha) \geq \mathbb{E}_\pi (1- \psi_\epsilon^\star) = \mathbb{E}_0 e^{\ln (L_\pi(Y))} (1- \psi_\epsilon^\star),$$
where $\psi_\epsilon^\star$ is the likelihood ratio test. In particular, $\psi_\epsilon^\star = 1$ if $\ln (L_\pi(Y))> t_\alpha$ where $t_\alpha$ is the $1-\alpha$ quantile of $\ln (L_\pi(Y))$ under $H_0$.
 Thanks to (\ref{eq:rapport_vrais}), we get 
$$ t_\alpha =-\frac{u_\epsilon^2(r_\epsilon)}{2} + u_\epsilon(r_\epsilon)H^{(\alpha)} + o(1) \ \mathrm{as} \ \epsilon \rightarrow 0,$$
which leads to the desired result, up to some simple algebra. Concerning the proof of (\ref{eq:rapport_vrais}), we refer to \cite{ISS_2012} for more details.  \\

Concerning the last part of the proof, we first have to compute the expectation and variance of the test statistics. 
Using a Markov inequality, it is then possible to prove that 
$$ \beta_\epsilon (\Theta(r_\epsilon),\Phi_{\epsilon,H^{(\alpha)}}) \rightarrow 0,$$
as soon as $u_\epsilon(r_\epsilon) \rightarrow +\infty$ as $\epsilon \rightarrow 0$. 
In the case where $u_\epsilon(r_\epsilon) = \mathcal{O}(1)$ as $\epsilon \rightarrow 0$, the Central Limit Theorem (with Lyapunov's condition)  indicates that the test statistics $T$ is asymptotically Gaussian. A precise study of its associated variance and expectation leads to 
the desired result (see \cite{ISS_2012} for a similar study in the 1-dimensional case).

\subsection{Proof of Proposition \ref{prop:cas1}}
\label{s:preuve_cas1}

Let 
$$b^2_{l}= \prod_{j=1}^d  l_j^{-2t_j} \ \mathrm{and} \ a^2_{l}=\prod_{j=1}^d |l_j|^{2s_j}, \ \forall l\in \mathbb{N}_*^d,$$
Recall that
$$ J_1=2^d \sum_{l \in \N^d} l_1^{4t_1}\ldots  l_d^{4t_d}(1-A l_1^{2s_1}\ldots l_d^{2s_d})_+.$$
We set $A=R^{-2 \bar{s}} $ where $\bar{s}=s_1+\ldots + s_d$. We define for $u=(u_1,\ldots, u_d) \in (\R_{+*})^d$ and 
$ s=(s_1,\ldots, s_d) \in (\R_{+*})^d$,
$$ S(u,s,R)=  \sum_{l\in \N^d, \pa{\frac{l_1}{R}}^{s_1}\ldots\pa{\frac{l_d}{R}}^{s_d} \leq 1}
 l_1^{u_1}\ldots  l_d^{u_d}.$$
\begin{lemma} \label{lem1}
Let for $j=1,\ldots d$, $c_j=(1+u_j)/s_j$. \\
Assume that $c_1>c_2>\ldots >c_d$.  Then we have 
$$   S(u,s,R) \sim_{R\ri \infty} \frac{R^{\bar{s} c_1}}{1+u_1} \prod_{j=2}^d \zeta \pa{c_1s_j-u_j},$$
where $\zeta(s)=\sum_{l\geq 1} l^{-s}$ for all $s>1$.\\
\end{lemma}
\textsc{Proof of Lemma \ref{lem1}}. We will use the following inequalities  
\begin{equation} \label{enca}
\forall u\geq 0, \  0\leq \sum_{l=1}^R l^u -\int_0^R x^u dx \leq 2^u R^u.
\end{equation}
\begin{equation} \label{enca2}
\forall -1<u< 0, \  0\leq \sum_{l=1}^R l^u \leq \int_0^R x^u dx .
\end{equation}

Indeed, for $u\geq 0$, 

$$\int_0^R x^u dx \leq \sum_{l=1}^R l^u \leq \int_1^{R+1} x^u dx,$$
and,
\begin{eqnarray*}
  0\leq \sum_{l=1}^R l^u -\int_0^R x^u dx &\leq & \int_R^{R+1} x^u dx \\
&\leq& \frac{1}{u+1}\cro{(R+1)^{u+1}-R^{u+1}}\\
&\leq& 2^u R^u.
\end{eqnarray*}
To simplify the notations, we omit in the sums that $l_1,\ldots, l_d \in \N^*$. 
Now,
$$ S(u,s,R)=  \sum_{ l_2^{s_2}\ldots l_d^{s_d} \leq R^{\bar{s}}}
 l_2^{u_2}\ldots  l_d^{u_d}\sum_{l_1=1}^{R^{\frac{\bar{s}}{s_1}} l_2^{-\frac{s_2}{s_1}}\ldots l_d^{-\frac{s_d}{s_1}}}
l_1^{u_1}  .$$ 
From Equation (\ref{enca}), we have 
\begin{eqnarray*}
 0\leq && S(u,s,R) -  \sum_{ l_2^{s_2}\ldots l_d^{s_d} \leq R^{\bar{s}}}  l_2^{u_2}\ldots  l_d^{u_d} 
\pa{\frac{1}{u_1+1}} \pa{R^{\frac{\bar{s}}{s_1}} l_2^{-\frac{s_2}{s_1}}\ldots l_d^{-\frac{s_d}{s_1}}}^{u_1+1}\\
&& \hspace{2cm}\leq 2^{u_1}  \sum_{ l_2^{s_2}\ldots l_d^{s_d} \leq R^{\bar{s}}}  l_2^{u_2}\ldots  l_d^{u_d} 
\pa{R^{\frac{\bar{s}}{s_1}} l_2^{-\frac{s_2}{s_1}}\ldots l_d^{-\frac{s_d}{s_1}}}^{u_1}.
\end{eqnarray*}
This leads to
\begin{eqnarray*}
 0\leq && S(u,s,R) -  \frac{R^{c_1\bar{s}}}{u_1+1}  \sum_{ l_2^{s_2}\ldots l_d^{s_d} \leq R^{\bar{s}} } 
l_2^{u_2-s_2c_1}\ldots l_d^{u_d-s_dc_1} \\
&& \hspace{2cm}\leq 2^{u_1} R^{\frac{\bar{s}}{s_1}u_1} \sum_{ l_2^{s_2}\ldots l_d^{s_d} \leq R^{\bar{s}}}  l_2^{u_2 -\frac{s_2}{s_1} u_1 }\ldots  
l_d^{u_d-\frac{s_d}{s_1} u_1}. 
\end{eqnarray*}
Note that for all $j \geq 2$, we have  $ u_j-s_jc_1 <-1$ since $c_j<c_1$. Hence the series 
$$ \sum_{ l_2^{s_2}\ldots l_d^{s_d} \leq R^{\bar{s}} } 
l_2^{u_2-s_2c_1}\ldots l_d^{u_d-s_dc_1} $$
converge. This leads to,
$$ \frac{R^{c_1\bar{s}}}{u_1+1}  \sum_{ l_2^{s_2}\ldots l_d^{s_d} \leq R^{\bar{s}} } 
l_2^{u_2-s_2c_1}\ldots l_d^{u_d-s_dc_1} \sim_{R\ri \infty} \frac{R^{\bar{s} c_1}}{1+u_1}
 \prod_{j=2}^d \zeta \pa{c_1s_j-u_j}.$$
It remains to prove that
$$ R^{\frac{\bar{s}}{s_1}u_1} \sum_{ l_2^{s_2}\ldots l_d^{s_d} \leq R^{\bar{s}}}  l_2^{u_2 -\frac{s_2}{s_1} u_1 }\ldots  
l_d^{u_d-\frac{s_d}{s_1} u_1} =o\pa{R^{\bar{s} c_1} }, \ \mathrm{as} \ R\rightarrow +\infty,$$
which is equivalent to 
\begin{equation}
\sum_{ l_2^{s_2}\ldots l_d^{s_d} \leq R^{\bar{s}}}  l_2^{u_2 -\frac{s_2}{s_1} u_1 }\ldots  
l_d^{u_d-\frac{s_d}{s_1} u_1} =o\pa{R^{\bar{s} \pa{c_1-\frac{u_1}{s_1}}}}, \ \mathrm{as} \ R\rightarrow +\infty . 
\label{inter_lem1}
\end{equation}
Since this quantity is positive, we will only give an upper bound. More generally, let us consider the Property
 that we denote
by $\mathcal{P}_j$ :
$$  \Sigma_j= \sum_{ l_j^{s_j}\ldots l_d^{s_d} \leq R^{\bar{s}}}  l_j^{u_j-\frac{s_j}{s_{j-1}} u_{j-1} }\ldots  
l_d^{u_d-\frac{s_d}{s_{j-1}} u_{j-1}} =o\pa{R^{\bar{s} \pa{c_1-\frac{u_{j-1}}{s_{j-1}}}}} . $$
Note that since $c_1 > c_j$, we have that $ c_1-\frac{u_{j-1}}{s_{j-1}} >0$ . \\
We want to prove (\ref{inter_lem1}), namely that $\mathcal{P}_2$ holds. 
Let us first prove that $\mathcal{P}_d$ holds, namely that
$$\sum_{l_d^{s_d} \leq R^{\bar{s}}}  
l_d^{u_d-\frac{s_d}{s_{d-1}} u_{d-1}} =o\pa{R^{\bar{s} \pa{c_1-\frac{u_{d-1}}{s_{d-1}}}}} . $$
This result clearly holds if $v_d= u_d-\frac{s_d}{s_{d-1}} u_{d-1} \leq -1$. If $v_d
> -1$, we get from (\ref{enca}) and (\ref{enca2}) that
\begin{eqnarray*}
 0\leq \sum_{l_d^{s_d} \leq R^{\bar{s}}}  
l_d^{v_d} &\leq &\int_0^{R^{\bar{s}/s_d}}x^{v_d}dx + 2^{v_d} R^{\frac{\bar{s}}{s_d}v_d}\\
&\leq& R^{\frac{\bar{s}}{s_d}(1+v_d)}\pa{C+o(1)}\\
&\leq & R^{\bar{s}(c_d-\frac{u_{d-1}}{s_{d-1}})}\pa{C+o(1)}\\
&=& o\pa{R^{\bar{s} \pa{c_1-\frac{u_{d-1}}{s_{d-1}}}}}. 
\end{eqnarray*}
We now assume that $\mathcal{P}_{j+1}$ holds and we want to prove that  $\mathcal{P}_j$ holds.  We have
\begin{eqnarray*}
\Sigma_j & = &  \sum_{ l_j^{s_j}\ldots l_d^{s_d} \leq R^{\bar{s}}}  l_j^{u_j-\frac{s_j}{s_{j-1}} u_{j-1} }\ldots  
l_d^{u_d-\frac{s_d}{s_{j-1}} u_{j-1}} ,\\
& = & \sum_{ l_{j+1}^{s_{j+1}}\ldots l_d^{s_d} \leq R^{\bar{s}}} 
 l_{j+1}^{u_{j+1}-\frac{s_{j+1}}{s_{j-1}} u_{j-1} }\ldots  
l_d^{u_d-\frac{s_d}{s_{j-1}} u_{j-1}} \sum_{l_j=1}^{R^{\frac{\bar{s}}{s_j}} l_{j+1}^{-\frac{s_{j+1}}{s_j}}
\ldots l_{d}^{-\frac{s_{d}}{s_j}}}
l_j^{u_j-\frac{s_j}{s_{j-1}}u_{j-1}}. 
\end{eqnarray*}
The property $\mathcal{P}_j$ clearly holds if $u_j-\frac{s_j}{s_{j-1}}u_{j-1} \leq -1$ (in this case, 
for all $l \geq j$, ${u_l-\frac{s_l}{s_{j-1}} u_{j-1}} \leq -1$). If  $u_j-\frac{s_j}{s_{j-1}}u_{j-1} >-1$, 
we use (\ref{enca}) to obtain 
\begin{eqnarray*}
 \Sigma_j &_\leq &  \sum_{ l_{j+1}^{s_{j+1}}\ldots l_d^{s_d} \leq R^{\bar{s}}} 
 l_{j+1}^{u_{j+1}-\frac{s_{j+1}}{s_{j-1}} u_{j-1} }\ldots  
l_d^{u_d-\frac{s_d}{s_{j-1}} u_{j-1}} \int_0^{R^{\frac{\bar{s}}{s_j}} l_{j+1}^{-\frac{s_{j+1}}{s_j}}
\ldots l_{d}^{-\frac{s_{d}}{s_j}}} x^{u_j-\frac{s_j}{s_{j-1}}u_{j-1} } dx\\
& & +C  \sum_{ l_{j+1}^{s_{j+1}}\ldots l_d^{s_d} \leq R^{\bar{s}}} 
 l_{j+1}^{u_{j+1}-\frac{s_{j+1}}{s_{j-1}} u_{j-1} }\ldots  
l_d^{u_d-\frac{s_d}{s_{j-1}} u_{j-1}} \pa{R^{\frac{\bar{s}}{s_j}} l_{j+1}^{-\frac{s_{j+1}}{s_j}}
\ldots l_{d}^{-\frac{s_{d}}{s_j}}}^{u_j-\frac{s_j}{s_{j-1}}u_{j-1} },\\
&\leq&  \frac{1}{1+ u_j-\frac{s_j}{s_{j-1}}u_{j-1}} R^{ \bar{s}(c_j-\frac{u_{j-1}}{s_{j-1}})} 
\sum_{ l_{j+1}^{s_{j+1}}\ldots l_d^{s_d} \leq R^{\bar{s}}} 
 l_{j+1}^{u_{j+1}- c_j s_{j+1}} \ldots  
l_d^{u_d-c_j {s_d}} \\
& & +C R^{\frac{\bar{s}}{s_j}(u_j-\frac{s_j}{s_{j-1}}u_{j-1})} \sum_{ l_{j+1}^{s_{j+1}}\ldots l_d^{s_d} 
\leq R^{\bar{s}}} 
 l_{j+1}^{u_{j+1}- u_j \frac{s_{j+1}}{s_j}}  \ldots  
l_d^{u_d- u_j \frac{s_d}{s_j}}.
\end{eqnarray*}
Since $ u_{j+1}- c_j s_{j+1} <-1$, the first series converges as $R\rightarrow +\infty$ and the first term is $O(R^{ \bar{s}(c_j-\frac{u_{j-1}}{s_{j-1}})})$
which is $o(R^{ \bar{s}(c_1-\frac{u_{j-1}}{s_{j-1}})})$ since $ c_j<c_1$. The second sum is $\Sigma_{j+1} $ which is
$o( R^{ \bar{s}(c_1-\frac{u_{j}}{s_{j}})})$ as $R\rightarrow +\infty$ since we have assumed that $\mathcal{P}_{j+1}$ holds. We then obtain that the second term is 
$o(R^{ \bar{s}(c_1-\frac{u_{j-1}}{s_{j-1}})})$ as $R\rightarrow +\infty$, which leads to the property  $\mathcal{P}_j$. 
\begin{flushright}
$\Box$
\end{flushright}

\begin{lemma} \label{J1}
 Let 
$$ J_1(R) =2^d \sum_{l_1,\ldots, l_d \in \N^*} l_1^{4t_1}\ldots  l_d^{4t_d}(1-R^{-2 \bar{s}} l_1^{2s_1}\ldots
 l_d^{2s_d})_+.$$
where $\bar{s}=s_1+\ldots + s_d$.
Let for $j=1,\ldots d$, $c_j=(1+4t_j)/s_j$. \\
Assume that $c_1>c_2>\ldots >c_d$.  Then we have 
$$  J_1(R) \sim_{R\ri \infty}2^d R^{\bar{s} c_1} \pa{\frac{2s_1}{(1+4t_1)(1+4t_1+2s_1)} }
\prod_{j=2}^d \zeta \pa{ c_1 s_j-4t _j}.$$
\end{lemma}
\textsc{Proof of Lemma \ref{J1}}. The proof follows directly from Lemma \ref{lem1} and easy computations by noticing that
$$  J_1(R) =2^d \pa{S(4t,s,R)- R^{-2 \bar{s}} S(4t+2s,s,R)}.$$
\begin{flushright}
$\Box$
\end{flushright}

\begin{lemma} \label{J2}
 Let 
$$ J_2(R) =2^d R^{-2 \bar{s}}  \sum_{l_1,\ldots, l_d \in \N^*} l_1^{4t_1+2s_1}\ldots  l_d^{4t_d+2s_d}   (1-R^{-2 \bar{s}} l_1^{2s_1}\ldots
 l_d^{2s_d})_+.$$
where $\bar{s}=s_1+\ldots + s_d$.
Let for $j=1,\ldots d$, $c_j=(1+4t_j)/s_j$. \\
Assume that $c_1>c_2>\ldots >c_d$.  Then we have 
$$  J_2(R) \sim_{R\ri \infty}2^d R^{\bar{s} c_1} \pa{\frac{2s_1}{(1+4t_1+2s_1)(1+4t_1+4s_1)} }
\prod_{j=2}^d \zeta \pa{ c_1 s_j-4t _j}.$$
\end{lemma}
\textsc{Proof of Lemma \ref{J2}}. The proof follows directly from Lemma \ref{lem1} and easy computations  by noticing that
$$  J_2(R) =2^d  R^{-2 \bar{s}} \pa{S(4t+2s,s,R)- R^{-2 \bar{s}} S(4t+4s,s,R)}.$$
\begin{flushright}
$\Box$
\end{flushright}

\begin{lemma} \label{J0}
 Let 
$$ J_0(R) =2^d \sum_{l_1,\ldots, l_d \in \N^*} l_1^{4t_1}\ldots  l_d^{4t_d}(1-R^{-2 \bar{s}} l_1^{2s_1}\ldots
 l_d^{2s_d})_+^2=J_1(R)-J_2(R).$$
where $\bar{s}=s_1+\ldots + s_d$.
Let for $j=1,\ldots d$, $c_j=(1+4t_j)/s_j$. \\
Assume that $c_1>c_2>\ldots >c_d$.  Then we have 
$$  J_0(R) \sim_{R\ri \infty}2^d R^{\bar{s} c_1} \pa{\frac{8s_1^2}{(1+4t_1)(1+4t_1+2s_1)(1+4t_1+4s_1)} }
\prod_{j=2}^d \zeta \pa{ c_1 s_j-4t _j}.$$
\end{lemma}

From the above lemmas, we can derive a separation rate for signal detection in this framework. Indeed,
$$r_{\epsilon}^2 = R^{-2 \bar{s}} \frac{J_1(R)}{J_2(R)}=  R^{-2 \bar{s}} \frac{D_1}{D_2},$$
where $  J_i(R) \sim_{R\ri \infty} D_i R^{\bar{s} c_1} $ for $i=0,1,2$. Then
$$ u_{\epsilon}^2(r_{\epsilon}) = \pa{\frac{r_{\epsilon}^4}{\epsilon^4} }\frac{J_0(R)}{2J_1^2(R)} = C
 \frac{r_{\epsilon}^{4+c_1}}{\epsilon^4} $$
where $C= \frac{D_0}{2D_1^2} \pa{ \frac{D_2}{D_1} }^{c_1/2}$. In particular 
$  u_{\epsilon}^2 (r_{\epsilon})  =\mathcal{O}( 1 ) $ for $r_{\epsilon}=r_{\epsilon}^* \sim \epsilon^{\frac{4}{4+c_1}}$ as $\epsilon \rightarrow 0$.

\subsection{Proof of Proposition \ref{prop:cas2}}
\label{s:preuve_cas2}

We consider the case where 
$$ b_l = \prod_{j=1}^d l_j^{-t_j} \ \mathrm{and} \ a_l = \prod_{j=1}^d e^{s l_j}, \ \forall l\in \N^d.$$
Remark that we are in an isotropic framework, i.e. $s_j = s$ for all $j\in \lbrace 1,\dots, d \rbrace$. As in the other cases, we start with the computation of $J_1$. Using our
 assumption and setting $A=e^{-2s u}$ we obtain
\begin{eqnarray*}
J_1 
& = & \sum_{l\in \N^d} b_l^{-4} (1-Aa_l^2),\\
& = & \sum_{l\in \N^d} \prod_{j=1}^d l_j^{4t_j} \left( 1 -  e^{2\sum_{j=1}^d s l_j -2s u} \right)_+,\\
& = & \sum_{\sum_{j=1}^d l_j \leq u} \prod_{j=1}^d l_j^{4t_j} \left( 1 - e^{2\sum_{j=1}^d s l_j -2s u} \right)_+,\\
& = & J_1' - J_1". 
\end{eqnarray*}
Simple algebra leads to
\begin{eqnarray*}
J_1' 
& := & \sum_{\sum_{j=1}^d l_j \leq u} \prod_{j=1}^d l_j^{4t_j},\\
& = & u^{\sum_{j=1}^d (4t_j +1)} \sum_{\sum_{j=1}^d (l_j/u) \leq 1} \prod_{j=1}^d \left(\frac{l_j}{u} \right)^{4t_j} \prod_{j=1}^d \frac{1}{u},\\
& = & C_1 u^{\sum_{j=1}^d (4t_j +1)} \int_{\mathbb{R}^d} \prod_{j=1}^d x_j^{4t_j} \mathbf{1}_{\lbrace \sum_{j=1}^d x_j \leq 1 \rbrace} dx
(1+o(1)) \ \mathrm{as} \ u \rightarrow +\infty.
\end{eqnarray*}
Now, we prove that $J_1"=J_1 \times o(1)$ as $u \rightarrow + \infty$. Let $\delta \in (0,1)$ a term whose value will be made precise later on. We can write that
\begin{eqnarray*}
J_1"
& := & \sum_{\sum_{j=1}^d l_j \leq u} \prod_{j=1}^d l_j^{4t_j} e^{2\sum_{j=1}^d s l_j -2s u},\\
& = & \sum_{\sum_{j=1}^d l_j \leq \delta u} \prod_{j=1}^d l_j^{4t_j} e^{2\sum_{j=1}^d s l_j -2s u} + \sum_{\delta u\leq \sum_{j=1}^d l_j \leq u} \prod_{j=1}^d l_j^{4t_j} e^{2\sum_{j=1}^d s l_j -2s u},\\
& \leq & e^{s(2\delta u -2u)}\sum_{\sum_{j=1}^d l_j \leq \delta u} \prod_{j=1}^d l_j^{4t_j} +  \sum_{\delta u\leq \sum_{j=1}^d l_j \leq u} \prod_{j=1}^d l_j^{4t_j},\\
& \sim & e^{-2s u(1-\delta)} \sum_{\sum_{j=1}^d l_j \leq \delta u} \prod_{j=1}^d l_j^{4t_j} + \prod_{j=1}^d u^{4t_j +1} 
\int_{\mathbb{R}^d} \prod_{j=1}^d x_j^{4t_j} \mathbf{1}_{\lbrace \delta \leq \sum_{j=1}^d x_j \leq 1 \rbrace} dx,\\
& = & J_1' \times o(1) \ \mathrm{as} \ u\rightarrow +\infty,
\end{eqnarray*}
setting for instance $\delta = \delta_u = 1-u^{-1/2}$.\\

The computation of $J_2$ follows essentially the same lines. First remark that
\begin{eqnarray*}
J_2
& = & \sum_{l\in \N^d} b_l^{-4} A a_l^2 (1-Aa_l^2)_+,\\
& = & \sum_{l\in \N^d} b_l^{-4} A a_l^2 - \sum_{l\in \N^d} b_l^{-4} A^2 a_l^4,\\
& := & J_2'+J_2".
\end{eqnarray*}
In a first time, we study $J_2'$:
\begin{eqnarray*}
J_2'
& := & \sum_{\sum_{j=1}^d l_j \leq u} \prod_{j=1}^d l_j^{4t_j} e^{2s ( \sum_{j=1}^d l_j - u)},\\
& = & \sum_{m=0}^u \sum_{\sum_{j=2}^d l_j \leq m} \left( m- \sum_{j=2}^d l_j \right)^{4t_1} \prod_{j=2}^d l_j^{4t_j} e^{2s (m-u)},\\
& = & \sum_{m=0}^u e^{-2s (u-m)} m^{\sum_{j=1}^d(4t_j +1) - 1} \sum_{\sum_{j=2}^d l_j/m \leq 1} \left( 1- \sum_{j=2}^d \frac{l_j}{m}
 \right)^{4t_1} \prod_{j=2}^d \left( \frac{l_j}{m} \right)^{4t_j} \frac{1}{m^{d-1}} ,\\
& = & \sum_{m=0}^u e^{-2s (u-m)} m^{\sum_{j=1}^d(4t_j +1) - 1} c_0(m),
\end{eqnarray*}
where
\begin{eqnarray*}
c_0(m) & = & \sum_{\sum_{j=2}^d l_j/m \leq 1} \left( 1- \sum_{j=2}^d \frac{l_j}{m} \right)^{4t_1} \prod_{j=2}^d
 \left( \frac{l_j}{m} \right)^{4t_j} \frac{1}{m^{d-1}},\\ 
& \sim & \int_{\mathbb{R_+}^d} \prod_{j=2}^d \left( x_j \right)^{4t_j} \pa{1-\sum_{j=1}^2 x_j}^{4t_1} \mathbf{1}_{\lbrace 
\sum_{j=2}^d x_j \leq 1 \rbrace}dx, \ \mathrm{as} \ m \rightarrow +\infty.
\end{eqnarray*}
Then, setting $l=u-m$, we obtain
\begin{eqnarray*}
J_2'
& = & \sum_{l=0}^u e^{-2s l} (u-l)^{\sum_{j=1}^d(4t_j +1) - 1} c_0(u-l),\\
& = & u^{\sum_{j=1}^d(4t_j +1) - 1} \sum_{l=0}^u e^{-2s l} \left(1-\frac{l}{u} \right)^{\sum_{j=1}^d(4t_j +1) - 1} c_0(u-l).
\end{eqnarray*}
It is possible to prove that the sum in the above formula converges as $u\rightarrow +\infty$. Using the same kind of algebra, 
we can also prove that $J_2"=J_2'\times o(1)$ as $u\rightarrow +\infty$. Therefore, we obtain the following asymptotic
$$ J_2 = C_2 u^{\sum_{j=1}^d(4t_j +1) - 1} (1+o(1)) \ \mathrm{as} \ u\rightarrow +\infty.$$
By the way, since $J_0 = J_1 - J_2$,
$$ J_0 = C_0 u^{\sum_{j=1}^d(4t_j +1)}(1+o(1)) \ \mathrm{as} \ u\rightarrow +\infty.$$

Now, we have got all the required material in order to compute the separation rate associated to this problem. First remark that there exists a constant $C$
such that as $\epsilon\rightarrow 0$,
\begin{eqnarray*}
r_\epsilon^2 = A \frac{J_1}{J_2} 
& \Leftrightarrow & r_\epsilon^2 = \frac{C_1}{C_2} e^{-2s u} \frac{u^{\sum_{j=1}^d(4t_j +1)}}{u^{\sum_{j=1}^d(4t_j +1)-1}}(1+o(1)),\\
& \Leftrightarrow &  r_\epsilon^2 = \frac{C_1}{C_2} e^{-2s u} u(1+o(1)),\\
& \Leftrightarrow & \ln (r_\epsilon) =  \frac{1}{2} \ln (u \frac{C_1}{C_2} (1+o(1))) - s u,\\
\end{eqnarray*}
The solution of the above equation satisfies 
$u =\ln (r_\epsilon^{-1/s}) (1+o(1))$. 
Then, as $\epsilon $ tends to $0$, 
\begin{eqnarray*}
u_\epsilon^2(r_\epsilon) = \left( \frac{r_\epsilon}{\epsilon} \right)^4 \frac{J_0}{2J_1^2} = O(1)
& \Leftrightarrow & \left( \frac{r_\epsilon}{\epsilon} \right)^4 u^{-\sum_{j=1}^d(4t_j +1)} = O(1),\\
& \Leftrightarrow & \left( \frac{r_\epsilon}{\epsilon} \right)^4 \left( \ln(r_\epsilon^{-1}) \right)^{-\sum_{j=1}^d(4t_j +1)} = O(1),,\\
& \Leftrightarrow & r_\epsilon \sim \epsilon \left( \ln(\epsilon^{-1}) \right)^{\sum_{j=1}^d(t_j +1/4)}.
\end{eqnarray*}

\subsection{Proof of Proposition \ref{prop:cas3}}
\label{s:preuve_cas3}
In this case, we will assume that
$$ b_l = \prod_{j=1}^d e^{-t_j l_j} \ \mathrm{and} \ a_l = \prod_{j=1}^d e^{s_j l_j} \ \forall l\in \N^d.$$
Moreover, we suppose that
$$ \frac{t_1}{s_1} > \frac{t_j}{s_j}, \ \forall j\in \lbrace 2,\dots,d \rbrace.$$
We start with the computation of the term $J_1$ defined as
\begin{eqnarray*}
J_1 & = & \sum_{l\in \N^d} b_l^{-4} (1- A a_l^2)_+,\\
& = & \sum_{l\in \N^d} e^{4 \sum_{j=1}^d t_j l_j} \left( 1- A e^{2\sum_{j=1}^d s_j l_j} \right)_+,\\
& = & \sum_{\sum_{j=1}^d s_j l_j \leq u} e^{4 \sum_{j=1}^d t_j l_j} \left( 1- A e^{2\sum_{j=1}^d s_j l_j} \right),
\end{eqnarray*}
setting $A = e^{-2u}$. Then
$$ J_1 = \sum_{\sum_{j=1}^d s_j l_j \leq u} e^{4 \sum_{j=1}^d t_j l_j} - A \sum_{\sum_{j=1}^d s_j l_j \leq u} e^{ \sum_{j=1}^d l_j(4t_j+2s_j)} := J_1' + J_1".$$
In a first time, we study $J_1'$. Remark that
$$ \sum_{j=1}^d s_j l_j \leq u  \Leftrightarrow l_1 \leq \frac{1}{s_1} \left( u - \sum_{j=2}^d s_j l_j \right).$$
Hence
\begin{eqnarray*}
J_1' & := & \sum_{\sum_{j=1}^d s_j l_j \leq u} e^{4 \sum_{j=1}^d t_j l_j},\\ 
& = & \sum_{\sum_{j=2}^d s_j l_j \leq u} e^{4\sum_{j=2}^d t_j l_j} \sum_{l_1=0}^{s_1^{-1}(u-\sum_{j=2}^d s_j l_j)} e^{4t_1 l_1},\\
& = & \sum_{\sum_{j=2}^d s_j l_j \leq u} e^{4\sum_{j=2}^d t_j l_j} e^{4\frac{t_1}{s_1} (u - \sum_{j=2}^d s_j l_j)} \sum_{l_1=0}^{s_1^{-1}(u-\sum_{j=2}^d s_j l_j)} e^{4t_1(l_1 -s_1^{-1}(u-\sum_{j=2}^d s_j l_j))} ,\\
& = & \sum_{\sum_{j=2}^d s_j l_j \leq u} e^{4\sum_{j=2}^d t_j l_j} e^{4\frac{t_1}{s_1} (u - \sum_{j=2}^d s_j l_j)} \mathcal{H}\left( (u-\sum_{j=2}^d s_j l_j),s_1,t_1\right), 
\end{eqnarray*}
where
$$ \mathcal{H}\left( x,s_1,t_1\right):= \sum_{l_1=0}^{s_1^{-1}x} e^{4t_1(l_1 -s_1^{-1}x)}= \mathcal{O}(1) \ \mathrm{as} \ x\rightarrow +\infty.$$
We obtain
\begin{eqnarray*}
J_1' & = & e^{4\frac{t_1}{s_1} u} \sum_{\sum_{j=2}^d s_j l_j \leq u} e^{4\sum_{j=2}^d t_j l_j} e^{-4\frac{t_1}{s_1}  \sum_{j=2}^d s_j l_j)} \mathcal{H}\left( (u-\sum_{j=2}^d s_j l_j),s_1,t_1\right),\\
& = & e^{4\frac{t_1}{s_1} u} \sum_{\sum_{j=2}^d s_j l_j \leq u} e^{-4\sum_{j=2}^d s_j l_j \left(\frac{t_1}{s_1} - \frac{t_j}{s_j} \right)} \mathcal{H}\left( (u-\sum_{j=2}^d s_j l_j),s_1,t_1\right),\\
& = & e^{4\frac{t_1}{s_1} u} c_1(u),
\end{eqnarray*}
where $c_1(u)=\mathcal{O}(1)$ as $u\rightarrow +\infty$. We are now interested in the computation of the term $J_1"$ defined as
\begin{eqnarray*}
J_1"
& = & A \sum_{\sum_{j=1}^d s_j l_j \leq u} e^{ \sum_{j=1}^d l_j(4t_j+2s_j)},\\
& = & A \sum_{\sum_{j=2}^d s_j l_j \leq u} e^{ \sum_{j=2}^d l_j(4t_j+2s_j)} \sum_{l_1=0}^{s_1^{-1}(u-\sum_{j=2}^d s_j l_j)} e^{(4t_1 + 2s_1)l_1},\\
& = & A \sum_{\sum_{j=2}^d s_j l_j \leq u} e^{ \sum_{j=2}^d l_j(4t_j+2s_j)} e^{(4t_1 + 2s_1)s_1^{-1}(u-\sum_{j=2}^d s_j l_j)}\\
& & \hspace{2cm}\times \sum_{l_1=0}^{s_1^{-1}(u-\sum_{j=2}^d s_j l_j)} e^{(4t_1 + 2s_1)(l_1-s_1^{-1}(u-\sum_{j=2}^d s_j l_j))},\\
& = & A e^{4\frac{t_1}{s_1} u}e^{2u} \sum_{\sum_{j=2}^d s_j l_j \leq u} e^{ \sum_{j=2}^d l_j(4t_j+2s_j)} e^{-(4t_1 + 2s_1)s_1^{-1}\sum_{j=2}^d s_j l_j}\\
& & \hspace{2cm} \times \mathcal{H}_1\left( (u-\sum_{j=2}^d s_j l_j),s_1,t_1 \right),
\end{eqnarray*}
where
$$ \mathcal{H}_1( x,s_1,t_1) =  \sum_{l_1=0}^{x} e^{(4t_1 + 2s_1)(l_1-x)}=\mathcal{O}(1), \ \mathrm{as} \ x\rightarrow +\infty.$$
Therefore, since $A=e^{-2u}$
\begin{eqnarray*}
J_1" & = & e^{4\frac{t_1}{s_1} u} \sum_{\sum_{j=2}^d s_j l_j \leq u} e^{ \sum_{j=2}^d l_j\lbrace (4t_j+2s_j)-(4t_1 + 2s_1)s_1^{-1}s_j \rbrace}  \mathcal{H}_1\left( (u-\sum_{j=2}^d s_j l_j),s_1,t_1 \right),\\
& = & e^{4\frac{t_1}{s_1} u}\sum_{\sum_{j=2}^d s_j l_j \leq u} e^{-4\sum_{j=2}^d s_j l_j \left(\frac{t_1}{s_1} - \frac{t_j}{s_j} \right)} \mathcal{H}_1\left( (u-\sum_{j=2}^d s_j l_j),s_1,t_1 \right),\\
& = & e^{4\frac{t_1}{s_1} u} c_1(u),
\end{eqnarray*}
where $c_1(u)=\mathcal{O}(u)$ as $u\rightarrow +\infty$. We finally obtain the asymptotic
$$ J_1 \sim e^{4\frac{t_1}{s_1} u}, \ \mathrm{as} \ u\rightarrow +\infty.$$
Using the same algebra, we obtain 
$$ J_2 \sim J_0 \sim e^{4\frac{t_1}{s_1} u}, \ \mathrm{as} \ u\rightarrow +\infty.$$
We can now study the separation rate associated to this framework. Since $J_1$ and $J_2$ are of the same order, we get that $A = r_\epsilon^2$. Then, as $\epsilon \rightarrow 0$,
\begin{eqnarray*}
u_\epsilon^2(r_\epsilon) =O(1)
& \Leftrightarrow & \left( \frac{r_\epsilon}{\epsilon} \right)^4 e^{-2\frac{t_1}{s_1}u} =O(1),\\
& \Leftrightarrow & \left( \frac{r_\epsilon}{\epsilon} \right)^4 A^{2\frac{t_1}{s_1}u} =O(1),\\
& \Leftrightarrow & r_\epsilon^{4(1+\frac{t_1}{s_1})} \sim \epsilon^4,\\
& \Leftrightarrow & r_\epsilon \sim \epsilon^{\frac{s_1}{s_1+t_1}}.
\end{eqnarray*}

\subsection{Proof of Proposition \ref{prop:cas4}}
\label{s:preuve_cas4}

In order to establish the separation rates related to this setting, we first need the following lemma.

\begin{lemma}
\label{lem:debut_cas4}
Let $d \in \mathbb{N}$ be fixed and assume that 
$$ t_1 > t_j \quad \forall j\in \lbrace 2,\dots, d \rbrace.$$
Then, there exists a constant $\mathcal{C}_d$ such that
\begin{equation}
 \sum_{\prod_{j=1}^d l_j \leq S} e^{4\sum_{j=1}^d t_j l_j} = \mathcal{C}_d e^{4t_1 S} (1+o(1)) \ \mathrm{as} \ S\rightarrow +\infty.
\label{eq:prop_cas4}
\end{equation}
\end{lemma}
\textsc{Proof}. Clearly, (\ref{eq:prop_cas4}) holds when $d$ is equal to $1$. Now, assume that equation (\ref{eq:prop_cas4}) holds for $d-1$. In such case, we prove that the same property holds for $d$, which will complete the proof of the lemma. In a first time, remark that
\begin{eqnarray*}
\sum_{\prod_{j=1}^d l_j \leq S}e^{4\sum_{j=1}^d t_j l_j}
& = & \sum_{\prod_{j=1}^{d-1} l_j \leq S}e^{4\sum_{j=1}^{d-1} t_j l_j} \times \sum_{l_d=1}^{S\prod_{j=1}^{d-1} l_j^{-1}} e^{4t_d l_d},\\
& = & \sum_{\prod_{j=1}^{d-1} l_j \leq S}e^{4\sum_{j=1}^{d-1} t_j l_j+ 4t_d S\prod_{j=1}^{d-1} l_j^{-1}} \times \sum_{l_d=1}^{S\prod_{j=1}^{d-1} l_j^{-1}} e^{4t_d (l_d - S\prod_{j=1}^{d-1} l_j^{-1})},\\
& = & \sum_{\prod_{j=1}^{d-1} l_j \leq S}e^{4\sum_{j=1}^{d-1} t_j l_j+ 4t_d S\prod_{j=1}^{d-1} l_j^{-1}} \times c_1\left(S\prod_{j=1}^{d-1} l_j^{-1} \right), 
\end{eqnarray*}
where
$$ c_1(U) := \sum_{l_d=1}^{U} e^{4t_d (l_d - U) } = \mathcal{O}(1) \ \mathrm{as} \ U\rightarrow +\infty.$$
Then, given a constant $\gamma \in ]0,1[$, the sum of interest can be decomposed as follows
\begin{eqnarray*}
\sum_{\prod_{j=1}^d l_j \leq S}e^{4\sum_{j=1}^d t_j l_j}
& = & \sum_{1\leq \prod_{j=1}^{d-1} l_j \leq (1-\gamma)S}e^{4\sum_{j=1}^{d-1} t_j l_j+ 4t_d S\prod_{j=1}^{d-1} l_j^{-1}} \times c_1\left(S\prod_{j=1}^{d-1} l_j^{-1} \right)\\
& & + \sum_{(1-\gamma)S \leq \prod_{j=1}^{d-1} l_j \leq S - \sqrt{S}}e^{4\sum_{j=1}^{d-1} t_j l_j+ 4t_d S\prod_{j=1}^{d-1} l_j^{-1}} \times c_1\left(S\prod_{j=1}^{d-1} l_j^{-1} \right)\\
& & \hspace{1cm} +
\sum_{S-\sqrt{S} \leq \prod_{j=1}^{d-1} l_j \leq S}e^{4\sum_{j=1}^{d-1} t_j l_j+ 4t_d S\prod_{j=1}^{d-1} l_j^{-1}} \times c_1\left(S\prod_{j=1}^{d-1} l_j^{-1} \right),\\
& := & T_1+T_2+T_3.
\end{eqnarray*}
We start we the control of the term $T_1$. Using simple algebra, we get that 
\begin{eqnarray*}
T_1 & := &  \sum_{1\leq \prod_{j=1}^{d-1} l_j \leq (1-\gamma)S}e^{4\sum_{j=1}^{d-1} t_j l_j+ 4t_d S\prod_{j=1}^{d-1} l_j^{-1}} \times c_1\left(S\prod_{j=1}^{d-1} l_j^{-1} \right),\\
& \leq & C e^{4t_d S} \sum_{1\leq \prod_{j=1}^{d-1} l_j \leq (1-\gamma)S}e^{4\sum_{j=1}^{d-1} t_j l_j} \leq C e^{4t_d S+ 4(1-\gamma)t_1S},
\end{eqnarray*}
where for the last inequality, we have used the hypothesis that equation (\ref{eq:prop_cas4}) holds for $d-1$.  Then, we can remark that
$$ T_1 = \mathcal{O}(e^{4t_d S+ 4(1-\gamma)t_1S}) = o(e^{4t_1 S}) \ \mathrm{as} \ S\rightarrow +\infty,$$
as soon as 
$$ 4t_1 (1-\gamma) S+4t_d S < 4t_1 S \Leftrightarrow  1>\gamma > \frac{t_d}{t_1}. $$
Now, we turn our attention to the term $T_2$. Assuming that equation (\ref{eq:prop_cas4}) holds for $d-1$, we get that
\begin{eqnarray*}
T_2 & := & 
\sum_{(1-\gamma)S \leq \prod_{j=1}^{d-1} l_j \leq S - \sqrt{S}}e^{4\sum_{j=1}^{d-1} t_j l_j+ 4t_d S\prod_{j=1}^{d-1} l_j^{-1}} \times c_1\left(S\prod_{j=1}^{d-1} l_j^{-1} \right),\\
& \leq & C \sum_{1 \leq \prod_{j=1}^{d-1} l_j \leq S - \sqrt{S}}e^{4\sum_{j=1}^{d-1} t_j l_j+ 4t_d (1-\gamma)^{-1}},\\
& = & \mathcal{O}(e^{4t_1(S-\sqrt{S})})=o(e^{4t_1 S}), \ \mathrm{as} \ S\rightarrow +\infty.
\end{eqnarray*}
Once again, if one assume that equation (\ref{eq:prop_cas4}) holds for $d-1$, we obtain that the term $T_3$ defined as 
\begin{eqnarray*}
T_3 & := & 
\sum_{S-\sqrt{S} \leq \prod_{j=1}^{d-1} l_j \leq S}e^{4\sum_{j=1}^{d-1} t_j l_j+ 4t_d S\prod_{j=1}^{d-1} l_j^{-1}} \times c_1\left(S\prod_{j=1}^{d-1} l_j^{-1} \right)
\end{eqnarray*}
can be surrunded as follows
\begin{eqnarray*}
& & e^{4t_d} \sum_{S-\sqrt{S} \leq \prod_{j=1}^{d-1} l_j \leq S}e^{4\sum_{j=1}^{d-1} t_j l_j} \times c_1\left(S\prod_{j=1}^{d-1} l_j^{-1} \right)\\
& & \hspace{2cm}  \leq T_3 \leq e^{4t_d S/(S-\sqrt{S})} \sum_{S-\sqrt{S} \leq \prod_{j=1}^{d-1} l_j \leq S}e^{4\sum_{j=1}^{d-1} t_j l_j} \times c_1\left(S\prod_{j=1}^{d-1} l_j^{-1} \right),\\
& \Leftrightarrow & T_3= \mathcal{C}_d e^{4t_1S}(1+o(1)) \ \mathrm{as} \ S\rightarrow +\infty.
\end{eqnarray*}
The proofs is a direct consequence on successive asymptotics of $T_1,T_2$ and $T_3$.
\begin{flushright}
$\Box$
\end{flushright}
\vspace{0.5cm}

Now, we can start the proof of Proposition \ref{prop:cas4} with the control of the term $J_1$. In a first time, set $S=A^{-1/2s}$ and remark that
\begin{eqnarray*}
J_1 & = & \sum_{l\in \mathbb{N}^d}  e^{4\sum_{j=1}^d t_j l_j} \left( 1-A \prod_{j=1}^d l_j^{2s} \right)_+,\\
& = & \sum_{\prod_{j=2}^{d} l_j \leq S} \sum_{l_1=1}^{S \prod_{j=2}^d l_j^{-1}}  e^{4\sum_{j=1}^d t_j l_j} \left( 1-A \prod_{j=1}^d l_j^{2s} \right)_+,\\
& = & \sum_{\prod_{j=2}^{d} l_j \leq S} e^{4\sum_{j=2}^{d} t_j l_j} \sum_{l_1=1}^{S \prod_{j=2}^d l_j^{-1}} e^{4t_1 l_1}   \left( 1-A \prod_{j=1}^d l_j^{2s} \right)_+,\\
& = & \sum_{\prod_{j=2}^{d} l_j \leq S} e^{4\sum_{j=2}^{d} t_j l_j + 4t_1  S\prod_{j=2}^d l_j^{-1}} \times c_0\left( S \prod_{j=2}^d l_j^{-1},t_1 \right),
\end{eqnarray*}
where
$$ c_0\left( R,t_1 \right) = \sum_{l=1}^{R} e^{4t_1( l-R)}   \left( 1- \left(\frac{l}{R} \right)^{2s} \right).$$
In particular, we have
\begin{eqnarray*}
c_0\left( R,t_1 \right)
& = & \sum_{m=1}^{R} e^{-4t_1 m}   \left( 1- \left(1-\frac{m}{R} \right)^{2s} \right),\\
& = & \sum_{m=1}^{\sqrt{R}} e^{-4t_1 m}   \left( 1- \left(1-\frac{m}{R} \right)^{2s} \right)+ \sum_{m=\sqrt{R}}^{R} e^{-4t_d m}   \left( 1- \left(1-\frac{m}{R} \right)^{2s} \right),\\
& = & \sum_{m=1}^{\sqrt{R}} e^{-4t_1 m}   \left( 2\frac{m}{R} + \mathcal{O}\left( \frac{m^2}{R^2} \right) \right) + \mathcal{O}\left( e^{-4t_1\sqrt{R}} \right), \ \mathrm{as} \ R \rightarrow +\infty,\\
& = & \mathcal{O}\left( R^{-1} \right), \ \mathrm{as} \ R \rightarrow +\infty.
\end{eqnarray*}
Now, we can concentrate our attention on the control of the term $J_1$. The sum in $J_1$ will be decomposed in two terms, namely
\begin{eqnarray*} 
J_1 & = & \left(  \sum_{1\leq \prod_{j=2}^{d} l_j \leq \sqrt{S}} + \sum_{\sqrt{S}\leq \prod_{j=2}^{d} l_j \leq S} \right) e^{4\sum_{j=2}^{d} t_j l_j + 4t_1  S\prod_{j=2}^d l_j^{-1}} \times c_0\left( S \prod_{j=2}^d l_j^{-1},t_1 \right),\\
& := & J_{1,1}+ J_{1,2} .
\end{eqnarray*}
Then, using Lemma \ref{lem:debut_cas4}, and the fact that the function $c_0$ is upper bounded,
\begin{eqnarray*} 
J_{1,2} 
& := & \sum_{\sqrt{S}\leq \prod_{j=2}^{d} l_j \leq S} e^{4\sum_{j=2}^{d} t_j l_j + 4t_1  S\prod_{j=1}^d l_j^{-1}}\times c_0\left( S \prod_{j=2}^d l_j^{-1},t_1 \right),\\
& \leq & C e^{4t_1 \sqrt{S}}\sum_{1\leq \prod_{j=2}^{d} l_j \leq S} e^{4\sum_{j=2}^{d} t_j l_j},\\
&  = & C e^{4t_1 \sqrt{S}} \times o(e^{4t_1 S}) = o(e^{4t_1 S})  \ \mathrm{as} \ S\rightarrow +\infty,
\end{eqnarray*}
Then,
\begin{eqnarray*}
J_{1,1}
&:=& \sum_{\prod_{j=2}^{d} l_j \leq \sqrt{S}} e^{4\sum_{j=2}^{d} t_j l_j + 4t_1  S\prod_{j=2}^d l_j^{-1}}\times c_0\left( S \prod_{j=2}^d l_j^{-1},t_1 \right),\\
& = & e^{4\sum_{j=2}^{d} t_j+4t_1S}  \times c_0\left( S \prod_{j=2}^d l_j^{-1},t_1 \right),\\
& & \hspace{3cm}         + \sum_{2 \leq \prod_{j=2}^{d} l_j \leq \sqrt{S}} e^{4\sum_{j=2}^{d} t_j l_j + 4t_1  S\prod_{j=2}^d l_j^{-1}}\times c_0\left( S \prod_{j=2}^d l_j^{-1},t_1 \right) ,\\
& = & J_{1,1}' +J_{1,1}".
%& \leq &  e^{4\sum_{j=2}^{d} t_j+4t_1S}\times \frac{1}{S}(1+o(1))  + \mathcal{O}\left( e^{2t_1 S}  \right)\times o(e^{4t_1S})  \ \mathrm{as} \ S\rightarrow +\infty.
\end{eqnarray*}
Then, we can remark that
$$ J_{1,1}' = \mathcal{C}_{1,d} \frac{e^{4t_1S}}{S}(1+o(1)), \ \mathrm{as} \ S \rightarrow +\infty,$$
while
$$ J_{1,1}" = \sum_{2 \leq \prod_{j=2}^{d} l_j \leq \sqrt{S}} e^{4\sum_{j=2}^{d} t_j l_j + 4t_1  S\prod_{j=2}^d l_j^{-1}}\times c_0\left( S \prod_{j=2}^d l_j^{-1},t_1 \right)  \leq C e^{2t_1S} \times o(e^{4t_1 \sqrt{S}}) = o(e^{4t_1 S}),$$
as $S\rightarrow +\infty$.
Finally,
$$ J_1 = \mathcal{C}_{1,d}\frac{e^{4t_1 S}}{S}(1+o(1)) \ \mathrm{as} \ S\rightarrow +\infty.$$

\vspace{1cm}
Using the same kind of algebra, we can then prove that
$$ J_2 := A\sum_{k\in \N^d} a_k^2 b_k^{-4} (1-Aa_k^2)_+ :=\mathcal{C}_{1,d}\frac{e^{4t_1 S}}{S}(1+o(1)) \ \mathrm{as} \ S\rightarrow +\infty,$$
and
$$ J_0 := \sum_{k\in \N^d}  b_k^{-4} (1-Aa_k^2)_+^2 :=\mathcal{C}_{0,d}\frac{e^{4t_1 S}}{S^2}(1+o(1)) \ \mathrm{as} \ S\rightarrow +\infty, $$
for some explicit constant $\mathcal{C}_{0,d}$.\\

In order to conclude the proof, we have to determine the minimax separation rate $r_\epsilon^\star$. Thanks to the asymptotics of $J_0,J_1$ and $J_2$ established above, we get 
$$ r_\epsilon^2 = A \frac{J_1}{J_2}= A(1+o(1)) \ \mathrm{as} \ A\rightarrow 0.$$ 
Hence, as $\epsilon \rightarrow 0$,
$$ u_\epsilon^2(r_\epsilon) \sim 1 \Leftrightarrow \epsilon^{-4} \frac{A^2}{J_2^2} J_0 \sim 1 \Leftrightarrow r_\epsilon^4 e^{-4t_1 r_\epsilon^{-1/2s}} \sim 1.$$
In particular,
$$ r_\epsilon^\star \sim \left( \frac{1}{4t_1} \ln \left( \frac{1}{\epsilon^4} \right)\right)^{-2s} \ \mathrm{as} \ \epsilon\rightarrow 0,$$
but we do not get sharp separation rates.

\subsection{Proof of Proposition \ref{prop:sobolev_cas1}}
\label{s:sobolev_cas1}
We here consider the case where
$$ b^2_{l}=\prod_{j=1}^d l_j^{-2t_j} \ \mathrm{and} \ a^2_{l}=\sum_{j=1}^d l_j^{2s_j}, \ \forall l \in \N^d.$$
We begin the proof with the study of $J_1$ defined as
$$ J_1=2^d \sum_{l_1,\ldots, l_d \in \N} \prod_{j=1}^d l_j^{4t_j} \left(1-A \sum_{j=1}^d l_j^{2s_j}\right)_+.$$
Setting $A=R_j^{-2 s_j} $ for all $j \in \lbrace 1,\ldots, d \rbrace$ and $\mathcal{D}=\prod_{j=1}^d R_j^{1+4t_j}$, we get
$$ J_1= 2^d \mathcal{D}  \sum_{l_1,\ldots, l_d \in \N^*} \prod_{j=1}^d \pa{\frac{l_j}{R_j}}^{4t_j}
 \pa{1-\sum_{j=1}^d \pa{\frac{l_j}{R_j}}^{2s_j} }_+  \prod_{j=1}^d \pa{\frac{1}{R_j}}.$$
Hence we have
$$  J_1\sim_{R_1,\ldots, R_d\ri \infty} \mathcal{D} 2^d  \int_{(\R^+)^d } \prod_{j=1}^d x_j^{4t_j} \left(1- \sum_{j=1}^d 
 x_j^{2s_j}\right)_+ dx_1\ldots dx_d. $$
We set $C_1= 2^d \int_{(\R^+)^d } \prod_{j=1}^d x_j^{4t_j} (1- \sum_{j=1}^d 
 x_j^{2s_j})_+dx_1\ldots dx_d. $
\begin{eqnarray*}
 C_1&=&  \frac{1}{ s_1\ldots s_d} \int_{(\R^+)^d } \prod_{j=1}^d v_j^{\frac{4t_j}{2s_j}}\left(1- \sum_{j=1}^d 
 v_j\right)_+  \prod_{j=1}^d v_j^{\frac{1}{2s_j}-1} dv_1\ldots dv_d \\
&=&  \frac{1}{ s_1\ldots s_d} \int_{T_d } \prod_{j=1}^d v_j^{\frac{4t_j+1}{2s_j}-1}\left(1- \sum_{j=1}^d 
 v_j \right)  dv_1\ldots dv_d 
\end{eqnarray*}
where $ T_d =\ac{(v_1,\ldots, v_d), v_j\geq 0, \sum_{j=1}^d v_j \leq 1}$. Let us recall Liouville's formula :
$$  \int_{T_d } \phi(v_1+\ldots + v_d)v_1^{p_1-1} \ldots v_d^{p_d-1} =\frac{\Gamma(p_1)\ldots \Gamma(p_d)}
{\Gamma(p_1+\ldots +p_d)} \int_0^1 \phi(u) u^{p_1+\ldots +p_d -1} du,$$
where $p_i >0$ for $i=1,\ldots,d$ and the integral in the right hand side is absolutely convergent. Using this
 formula, and setting $\tilde{s}= \sum_{j=1}^d (1+4t_j)/s_j$,
\begin{eqnarray*}
C_1&=&  \frac{ \prod_{j=1}^d  \Gamma\pa{\frac{4t_j+1}{2s_j} }}{ s_1\ldots s_d \Gamma(\tilde{s})}
\int_0^1 (1-u)u^{\tilde{s}-1} du \\
&=&  \frac{ \prod_{j=1}^d  \Gamma\pa{\frac{4t_j+1}{2s_j} }}{ s_1\ldots s_d \Gamma(\tilde{s})}
\frac{1}{\tilde{s}(\tilde{s}+1)}.
\end{eqnarray*}
Now, consider the term $J_2$ defined as
$$ J_2=2^d \sum_{l_1,\ldots, l_d \in \N^*} \prod_{j=1}^d l_j^{4t_j}  A \left( \sum_{j=1}^d l_j^{2s_j} \right) 
\left(1-A \sum_{j=1}^d l_j^{2s_j} \right)_+.$$
We 	set $A=R_j^{-2 s_j} $ for all $j=1,\ldots, d$. Setting $\mathcal{D}=\prod_{j=1}^d R_j^{1+4t_j}$, we have
$$ J_2= \mathcal{D}2^d   \sum_{l_1,\ldots, l_d \in \N^*} \prod_{j=1}^d \pa{\frac{l_j}{R_j}}^{4t_j}
 \left( \sum_{j=1}^d \pa{\frac{l_j}{R_j}}^{2s_j} \right)
 \pa{1-\sum_{j=1}^d \pa{\frac{l_j}{R_j}}^{2s_j} }_+  \prod_{j=1}^d \pa{\frac{1}{R_j}}.$$
Hence we have
$$  J_2\sim_{R_1,\ldots, R_d\ri \infty} \mathcal{D} 2^d  \int_{(\R^+)^d } \prod_{j=1}^d x_j^{4t_j} \left( \sum_{j=1}^d 
 x_j^{2s_j} \right) \left(1- \sum_{j=1}^d 
 x_j^{2s_j}\right)_+dx_1\ldots dx_d. $$
Setting 
$$C_2= 2^d \int_{(\R^+)^d } \prod_{j=1}^d x_j^{4t_j}( \sum_{j=1}^d 
 x_j^{2s_j})(1- \sum_{j=1}^d 
 x_j^{2s_j})_+dx_1\ldots dx_d. $$
we get, using similar computations as above
$$ C_2= \frac{ \prod_{j=1}^d  \Gamma\pa{\frac{4t_j+1}{2s_j} }}{ s_1\ldots s_d \Gamma(\tilde{s})}
\frac{1}{(\tilde{s}+1)(\tilde{s}+2)}.$$
In the same manner, we have 
$$  J_2\sim_{R_1,\ldots, R_d\ri \infty} \mathcal{D} C_0$$
with
$$ C_0=  \frac{ \prod_{j=1}^d  \Gamma\pa{\frac{4t_j+1}{2s_j} }}{ s_1\ldots s_d \Gamma(\tilde{s})}
\frac{2}{\tilde{s}(\tilde{s}+1)(\tilde{s}+2)}.$$
Let us now determine a separation rate in this framework. \\
$$r_{\epsilon}^2 = A \frac{J_1}{J_2}=  R^{-2 s_j} \frac{C_1}{C_2},$$
hence $ R_j= \pa{\frac{C_1}{C_2}}^{1/(2s_j)} r_{\epsilon}^{-1/s_j}$. 
$$ u_{\epsilon}^2(r_{\epsilon}) = \pa{\frac{r_{\epsilon}^4}{\epsilon^4} }\frac{J_0}{2J_1^2} = 
\pa{\frac{r_{\epsilon}^4}{\epsilon^4} }\frac{C_0}{2C_1^2 \mathcal{D}} .$$
$$\mathcal{D}=  \prod_{j=1}^d R_j^{1+4t_j}  =  \prod_{j=1}^d \pa{\frac{C_1}{C_2}}^{\frac{1+4t_j}{2s_j}}
r_{\epsilon}^{-\frac{(1+4t_j)}{s_j}} .$$
This leads to 
$$ u_{\epsilon}^2(r_{\epsilon})= C \epsilon^{-4}r_{\epsilon}^{ 4+ \sum_{j=1}^d \frac{(1+4t_j)}{s_j}}.$$
with
$C=\frac{C_0}{2C_1^2}\pa{\frac{C_2}{C_1}}^{\tilde{s}}.$
In particular
$  u_{\epsilon}^2(r_{\epsilon}) =\mathcal{O}( 1) $ for $r_{\epsilon}=r_{\epsilon}^*= 
\epsilon^{{4}/\pa{4+\sum_{j=1}^d \frac{(1+4t_j)}{s_j}}}.$

\subsection{Proof of Proposition \ref{prop:sobolev_cas2}}
\label{s:sobolev_cas2}
In this case, we have
$$ b_l = \prod_{j=1}^d e^{-t_j l_j} , \ \mathrm{and} \ a_l = \left(\sum_{j=1}^d l_j \right)^s, \ \forall l\in \N^d.$$
Remark that we consider an isotropic framework: the regularity is the same for all the $d$ directions.\\

We start with the computation of $J_1$ defined as
\begin{eqnarray*}
J_1  & = & \sum_{l\in \N^d} b_l^{-4} (1- A a_l^2 )_+,\\
& = & \sum_{l\in\N^d} e^{4\sum_{j=1}^d t_j l_j} \left( 1- A \left( \sum_{j=1}^d l_j \right)^{2s} \right)_+,\\
& = & \sum_{ \sum_{j=1}^d l_j \leq R} e^{4\sum_{j=1}^d t_j l_j} \left( 1- A \left( \sum_{j=1}^d l_j \right)^{2s} \right), 
\end{eqnarray*}
where $R:=A^{-1/2s}$. Setting $m=\sum_{j=1}^d l_j$, we get
\begin{eqnarray*}
J_1 & = & \sum_{m=0}^R \sum_{ \sum_{j=2}^d l_j \leq m} e^{4t_1 m + 4 \sum_{j=2}^d l_j (t_j - t_1)} (1-Am^{2s}),\\
& = & \sum_{m=0}^R e^{4t_1 m} (1-Am^{2s}) \left( \sum_{ \sum_{j=2}^d l_j \leq m} e^{4 \sum_{j=2}^d l_j (t_j - t_1)} \right),\\
& = & \sum_{m=0}^R e^{4t_1 m} (1-Am^{2s}) c_0(m),
\end{eqnarray*}
where
$$ c_0(m) = \sum_{ \sum_{j=2}^d l_j \leq m} e^{4 \sum_{j=2}^d l_j (t_j - t_1)}, \ \forall m\in \mathbb{N}.$$
Since we have assumed that $t_1 > t_j$ for all $j\geq 2$, $c_0(m)=\mathcal{O}(1)$ as $m\rightarrow +\infty$. Let $\delta \in [0,1]$ be a constant which will be made precise later on. We can write that
\begin{eqnarray*}
J_1 & = & \sum_{m=0}^R e^{4t_1 m} \left( 1-\left( \frac{m}{R} \right)^{2s} \right) c_0(m),\\
& = & e^{4t_1 R} \sum_{m=0}^R e^{4t_1 (m-R)} \left( 1-\left( 1+\frac{m-R}{R} \right)^{2s} \right) c_0(m),\\
& = & e^{4t_1 R} \sum_{l=0}^R e^{-4t_1 l} \left( 1-\left( 1 - \frac{l}{R} \right)^{2s} \right) \tilde c_0(l),\\
& = & e^{4t_1 R} \sum_{l=0}^{\delta R} e^{-4t_1 l} \left( 1-\left( 1 - \frac{l}{R} \right)^{2s} \right) \tilde c_0(l) + e^{4t_1 R} \sum_{l=\delta R}^R e^{-4t_1 l} \left( 1-\left( 1 - \frac{l}{R} \right)^{2s} \right) \tilde c_0(l),\\
& := & T_1 + T_2. 
\end{eqnarray*}
Concerning the term $T_2$, using simple algebra, we get that
\begin{eqnarray*}
T_2 & := &  e^{4t_1 R} \sum_{l=\delta R}^R e^{-4t_1 l} \left( 1-\left( 1 - \frac{l}{R} \right)^{2s} \right) \tilde c_0(l),\\
& \leq & e^{4t_1 R} C e^{-4t_1 \delta R} \leq C e^{4t_1 (1-\delta)R}.
\end{eqnarray*}
In order to compute $T_1$, we will use the Taylor expansion
$$ \left( 1- \frac{l}{R} \right)^{2s} = 1- 2s\frac{l}{R} + \mathcal{O} \left( \frac{l^2}{R^2} \right).$$
We obtain
\begin{eqnarray*}
T_1 & := & e^{4t_1 R} \sum_{l=0}^{\delta R} e^{-4t_1 l} \left( 1-\left( 1 - \frac{l}{R} \right)^{2s} \right) \tilde c_0(l),\\
& = & e^{4t_1 R} \sum_{l=0}^{\delta R} e^{-4t_1 l} \left( 2s\frac{l}{R} + \mathcal{O} \left( \frac{l^2}{R^2} \right) \right) \tilde c_0(l),\\
& = & \frac{2s}{R} e^{4t_1 R} \left( \sum_{l=0}^{\delta R} e^{-4t_1 l}l \tilde c_0(l) + \frac{\mathcal{O}(1)}{2sR} \sum_{l=0}^{\delta R} e^{-4t_1 l}l^2 \tilde c_0(l) \right),\\
& = &  \frac{2s}{R} e^{4t_1 R} \Delta(R,t_1)(1+o(1)).
\end{eqnarray*}
Hence 
$$ J_1 = 2s \Delta(R,t_1) \frac{e^{4t_1 R}}{R} (1+o(1)), \ \mathrm{as} \ R\rightarrow +\infty,$$
since $\Delta(R,t_1)= \mathcal{O}(1)$ as $R\rightarrow +\infty$.
Using the same algebra, we obtain
\begin{eqnarray*}
J_2 & = & \sum_{l\in \N^d} b_l^{-4} Aa_l^2(1- A a_l^2 )_+,\\
& = & \sum_{m=0}^R e^{4mt_1} \left( \frac{m}{R} \right)^{2s} \left( 1- \left(\frac{m}{R} \right)^{2s} \right) c_0(m),\\
& = & \sum_{m=0}^R e^{4mt_1} \left( \frac{m}{R} \right)^{2s} c_0(m) - \sum_{m=0}^R e^{4mt_1} \left( \frac{m}{R} \right)^{4s}c_0(m),\\
& = & \sum_{m=0}^R e^{4mt_1} \left( 1- \left( \frac{m}{R} \right)^{2s} \right) c_0(m) - \sum_{m=0}^R e^{4mt_1} \left( 1-\left( \frac{m}{R} \right)^{4s} \right) c_0(m),\\
& = & 4s \Delta(R,t_1) \frac{e^{4t_1 R}}{R} (1+o(1)) - 2s \Delta(R,t_1) \frac{e^{4t_1 R}}{R} (1+o(1)),\\
& = & 2s \Delta(R,t_1) \frac{e^{4t_1 R}}{R} (1+o(1)).
\end{eqnarray*}
Remark that we obtain exactly the same asymptotics for $J_1$ and $J_2$ which indicates the presence of sharp separation rates. Concerning the term $J_0$, we can prove that
\begin{eqnarray*}
J_0 & = & \sum_{l\in \N^d} b_l^{-4} (1- A a_l^2 )^2_+,\\
& = & \sum_{m=0}^R e^{4t_1 m} \left( 1- \left(\frac{m}{R}\right) \right)^2 c_0(m),\\
& = & 4s^2 \tilde \Delta(R,t_1) \frac{e^{4t_1 R}}{R^2} (1+o(1)), \ \mathrm{as} \ R\rightarrow +\infty.
\end{eqnarray*}
In order to find the corresponding separation rates, we have to solve
$$ r_\epsilon^2 = A \frac{J_1}{J_2}, \ \mathrm{and} \ u_\epsilon^2(r_\epsilon) = \left( \frac{r_\epsilon}{\epsilon} \right)^4 \frac{J_0}{2J_1^2}.$$
First remark that
$$ r_\epsilon^2 = A \frac{J_1}{J_2} \Leftrightarrow r_\epsilon^2 \sim A \Leftrightarrow r_\epsilon^2 \sim R^{-2s} \Leftrightarrow R \sim r_\epsilon^{-1/s}.$$
Hence
\begin{equation*}
\left( \frac{r_\epsilon}{\epsilon} \right)^4 \frac{J_0}{2J_1^2} =  \frac{r_\epsilon^4}{\epsilon^4} e^{- 4t_1 r_{\epsilon}^{-1/s}} 
% \sim 1
%& \Leftrightarrow & r_\epsilon^4 \frac{1/R^2}{1/R^2} \frac{e^{4t_1 R}}{e^{8t_1 R}} \sim \epsilon^4,\\
%& \Leftrightarrow & r_\epsilon^4 e^{- 4t_1 r_{\epsilon}^{-1/s}} \sim \epsilon^4,\\
%& \Leftrightarrow & r_\epsilon^2 e^{- 2t_1 r_{\epsilon}^{-1/s}} \sim \epsilon^2,\\
%& \Leftrightarrow & r_\epsilon \sim  \left( \log \left(\frac{1}{\epsilon}\right) \right)^{-s}.
\end{equation*} 
Hence, if $ r_\epsilon = \pa{C\log(1/\epsilon)}^{-s} $ with $C>1/t_1$, then $ u_\epsilon^2 (r_\epsilon)\rightarrow +\infty$ as $\epsilon \rightarrow
0$.
If $ r_\epsilon = \pa{C\log(1/\epsilon)}^{-s} $ with $C\leq 1/t_1$, then $ u_\epsilon^2 (r_\epsilon)\rightarrow 0 $.

%\vspace{0.5cm}
%\noindent
%\underline{Remark:} Up to now, we have supposed that $t_1> t_j$ for all $j\geq 2$. What happends if there is equality? Assume for instance that $d=3$ and $\beta_2 = t_1 > \beta_3$. In such a case,
%$$ c_0(m) = \sum_{l_2+l_3 \leq m} e^{4l_3 (\beta_3 - t_1)} = \sum_{l_3=0}^m (m-l_3) e^{-4l_3(\beta_3 - t_1)} \sim m.$$
%This should be writen properly but it seems that in such a case
%$$ J_1 \sim J_2 \sim e^{4t_1 R}, \ \mathrm{and} \ J_0 \sim \frac{1}{R} e^{4t_1 R}.$$
%Hence
%\begin{eqnarray*}
%\left( \frac{r_\epsilon}{\epsilon} \right)^4 \frac{J_0}{2J_1^2} \sim 1
%& \Leftrightarrow & r_\epsilon^4 \frac{1}{R} e^{-4t_1 R} \sim \epsilon^4,\\
%& \Leftrightarrow & r_\epsilon^2 r_{\epsilon}^{1/2s} e^{- 2t_1 r_{\epsilon}^{-1/s}} \sim \epsilon^2,\\
%& \Leftrightarrow & r_\epsilon \sim  \left( \log \left(\frac{1}{\epsilon}\right) \right)^{-s}.
%\\end{eqnarray*} 
%This does not change anything regarding the rate...

\bibliography{Inverse}
\bibliographystyle{plain}

\end{document}